\documentclass{article}

\let\epsilon=\varepsilon{}
\usepackage[utf8]{inputenc}
\usepackage{mathtools}
\usepackage{algorithm}
\usepackage{myalgorithmic}
\usepackage{amsthm}
\usepackage{amssymb}
\usepackage{bm}
\usepackage{booktabs}
\usepackage{comment}
\usepackage{etoolbox}
\usepackage{float}
\usepackage{graphicx}
\usepackage{needspace}
\usepackage[group-separator={,}]{siunitx}
\usepackage{tikz}
\usepackage{multirow}
\usepackage{hyperref}
\usepackage{paralist}
\usepackage{xcolor}
\usepackage{authblk}
\usepackage{microtype}
\usepackage{placeins}
\usepackage{lmodern}
\usepackage[tableposition=top]{caption}
\usepackage[margin=1in]{geometry}
\usepackage[T1]{fontenc}

\allowdisplaybreaks  

\usetikzlibrary{3d,calc}

\usepackage{import}

\usepackage[backend=biber, style=numeric, giveninits=true]{biblatex}

\DeclareFieldFormat*{title}{\mkbibemph{#1}}
\DeclareFieldFormat*{citetitle}{\mkbibemph{#1}}
\DeclareFieldFormat{journaltitle}{#1}

\renewbibmacro*{in:}{%
  \ifentrytype{article}
    {}%
    {\printtext{\bibstring{in}\intitlepunct}}%
    }

\newbibmacro*{pubinstorg+location+date}[1]{%
  \printlist{#1}%
  \newunit%
  \printlist{location}%
  \newunit%
  \usebibmacro{date}%
  \newunit}

\renewcommand{\Re}{\mathop{\mathrm{Re}}%
}
\newtheorem{theorem}{Theorem}
\newtheorem{lemma}[theorem]{Lemma}

\newtheorem{definition}{Definition}

\newcommand{\ilist}[2]{%
  \ifstrequal{#1}{1}{U_{#2}}{%
  \ifstrequal{#1}{2}{V_{#2}}{%
  \ifstrequal{#1}{3}{W_{#2}}{%
  \ifstrequal{#1}{3close}{W^\mathrm{close}_{#2}}{%
  \ifstrequal{#1}{3far}{W^\mathrm{far}_{#2}}{%
  \ifstrequal{#1}{4}{X_{#2}}{%
  \ifstrequal{#1}{4close}{X^\mathrm{close}_{#2}}{%
  \ifstrequal{#1}{4far}{X^\mathrm{far}_{#2}}{}%
}}%
}}%
}}%
}}

\def\Mpole{\mathsf{M}}
\def\Locfar{\mathsf{L}^{\text{far}}%
}
\def\Potnear{\mathsf{P}^{\text{near}}%
}
\def\PotW#1{\mathsf{P}^{\ilist{3}{}}_{#1}}
\def\Lqbxnear#1{\mathsf{L}^{\text{qbx},\text{near}}_{#1}}
\def\LqbxW#1{\mathsf{L}^{\text{qbx},\ilist{3}{}}_{#1}}
\def\Lqbxfar#1{\mathsf{L}^{\text{qbx},\text{far}}_{#1}}
\def\Lqbx#1{\mathsf{L}^{\text{qbx}}_{#1}}

\def\ptpot{\Pi}

\newlength{\longrowlength}
\newlength{\figurewidth}
\newcommand{\nmax}{{n_{\mathrm{max}}%
}}

\newcommand{\pqbx}{{p_\mathrm{qbx}%
}}
\newcommand{\pfmm}{{p_\mathrm{fmm}%
}}
\newcommand{\pquad}{{p_\mathrm{quad}%
}}
\newcommand{\nmpole}{{n_\mathrm{mpole}%
}}

\newcommand{\closedbox}{\overline{B_\infty}}
\newcommand{\closedball}{\overline{B_2}}

\newcommand{\ancestors}{\mathsf{Ancestors}}
\newcommand{\descendants}{\mathsf{Descendants}}
\newcommand{\tcr}{\mathsf{TCR}}
\newcommand{\parent}{\mathsf{Parent}}
\newcommand{\adequatesep}{\prec}

\newcommand{\mapping}{\Psi}

\newcommand{\algbrand}{GIGAQBX}

\newcommand{\kernel}{\mathcal{K}}

\newcommand{\norm}[1]{\lvert#1\rvert}
\newcommand{\norminf}[1]{\lvert#1\rvert_\infty}

\newcommand{\qbxsinglelayercont}[1]{(\mathcal{S}_\mathrm{global}^{(#1)})}

\newcommand{\localqbxsinglelayercont}[2]{(\mathcal{S}_{\mathrm{local},#2}^{(#1)})}

\newcommand{\qbxsinglelayer}[1]{(\tilde{\mathcal{S}}_\mathrm{global}^{(#1)})}

\makeatletter
\tikzoption{canvas is plane}[]{\@setOxy#1}
\def\@setOxy O(#1,#2,#3)x(#4,#5,#6)y(#7,#8,#9)%
  {\def\tikz@plane@origin{\pgfpointxyz{#1}{#2}{#3}}%
   \def\tikz@plane@x{\pgfpointxyz{#4}{#5}{#6}}%
   \def\tikz@plane@y{\pgfpointxyz{#7}{#8}{#9}}%
   \tikz@canvas@is@plane}
\makeatother

\newcommand\DrawSphere[5][black]{%
\begin{scope}[canvas is plane={O(#2,#3,#4)x(#2+#5,#3,#4)y(#2,#3+#5,#4)}]
  \draw[#1] (0,0) circle (1);
\end{scope}
\begin{scope}[canvas is plane={O(#2,#3,#4)x(#2,#3,#4+#5)y(#2,#3+#5,#4)}]
  \draw[#1] (0,-1) arc (270:270+180:1);
  \draw[#1,dotted] (0,1) arc (90:270:1);
\end{scope}
\begin{scope}[canvas is plane={O(#2,#3,#4)x(#2+#5,#3,#4)y(#2,#3,#4+#5)}]
  \draw[#1] (1,0) arc (0:180:1);
  \draw[#1,dotted] (-1,0) arc (180:360:1);
\end{scope}
\begin{scope}[canvas is plane={O(#2,#3,#4)x(#2+#5/sqrt 2,#3,#4+#5/sqrt 2)y(#2,#3+#5,#4)}]
  \draw[#1] (0,-1) arc (270:270+180:1);
  \draw[#1,dotted] (0,1) arc (90:270:1);
\end{scope}
\begin{scope}[canvas is plane={O(#2,#3,#4)x(#2-#5/sqrt 2,#3,#4+#5/sqrt 2)y(#2,#3+#5,#4)}]
  \draw[#1] (0,-1) arc (270:270+180:1);
  \draw[#1,dotted] (0,1) arc (90:270:1);
\end{scope}
}

\tikzset{%
  >=latex, 
  coord/.style={inner sep=0pt, outer sep=0pt, minimum size=3pt, circle},
}

\definecolor{qbxcolor}{RGB}{43,131,186}
\definecolor{srccolor}{RGB}{63,140,52}
\definecolor{localcolor}{RGB}{215,25,28}
\colorlet{cubecolor}{gray}

\newenvironment{algbreakable}[1]{
    \refstepcounter{algorithm}
    \needspace{3\baselineskip}
    \noindent\rule{\textwidth}{0.5pt}
    \textbf{Algorithm~\arabic{algorithm}:} #1\\
    \rule[1.25ex]{\textwidth}{0.5pt}
    \vspace{-4ex}
  }
  {
    \vspace{1ex}
    \rule{\textwidth}{0.5pt}
  }

\def\algstage#1#2{\vspace{2ex}\begin{minipage}{0.97\textwidth}\STATE{\textit{#1}}#2\end{minipage}}

\newcommand{\triforce}{\begin{tikzpicture}
  [scale=0.15,y={(0.5cm,0.866cm)}] 
  \draw ++(0,0)--++(1,0)--++(-1,1)--cycle;
  \draw ++(1,0)--++(1,0)--++(-1,1)--cycle;
  \draw ++(0,1)--++(1,0)--++(-1,1)--cycle;
\end{tikzpicture}%
}

\def\sqrtthree{1.732}

\newcommand{\ditto}[1][.4pt]{---~\textquotedbl~---}

\begin{document}

\title{Optimization of Fast~Algorithms for Global~Quadrature~by~Expansion
  Using Target-Specific~Expansions}

\author[1]{Matt Wala\thanks{\texttt{wala1@illinois.edu}}}
\author[1]{Andreas Klöckner\thanks{\texttt{andreask@illinois.edu}}}
\affil[1]{Department of Computer Science, University of Illinois at Urbana-Champaign}

\def\today{November 14, 2019}
\date{\today}

\maketitle

\begin{abstract}
  We develop an algorithm for the asymptotically fast evaluation of layer
  potentials close to and on the source geometry, combining Geometric Global
  Accelerated QBX (`GIGAQBX') and target-specific expansions.
  GIGAQBX is a fast high-order scheme for evaluation of layer potentials based
  on Quadrature by Expansion (`QBX') using local expansions formed via the Fast
  Multipole Method (FMM).
  Target-specific expansions serve to lower the cost of the formation and evaluation of QBX local expansions, reducing the
  associated computational effort from $O((p+1)^{2})$ to $O(p+1)$ in
  three dimensions, without any accuracy loss compared with conventional
  expansions, but with the loss of source/target separation in the expansion
  coefficients. GIGAQBX is a `global' QBX scheme, meaning that the potential is
  mediated entirely through expansions for points close to or on the
  boundary. In our scheme, this single global expansion is
  decomposed into two parts that are evaluated separately: one part incorporating
  near-field contributions using target-specific expansions, and one
  part using conventional spherical harmonic expansions of far-field contributions,
  noting that convergence guarantees only exist for the sum of the two sub-expansions.
  By contrast, target-specific expansions were originally introduced as an
  acceleration mechanism for `local' QBX schemes, in which the far-field
  does not contribute to the QBX expansion. Compared with the
  unmodified GIGAQBX algorithm, we show through a reproducible, time-calibrated cost model
  that the combined scheme yields a considerable cost reduction for the near-field
  evaluation part of the computation. We support the effectiveness of
  our scheme through numerical results demonstrating performance improvements
  for Laplace and Helmholtz kernels.
\end{abstract}


\newenvironment{DIFnomarkup}{}{\ignorespacesafterend}

\section{Introduction}%
\label{sec:intro}

The numerical realization of integral equation methods for the
solution of boundary value problems of elliptic partial differential equations
presents a number of technical challenges. Chief among them is the
\emph{accurate} and \emph{rapid} evaluation of \emph{layer potentials},
such as the single-layer potential
\begin{equation}
  \label{eqn:slp}
  (\mathcal{S} \mu)(x) \coloneqq \int_\Gamma \kernel(x,y) \mu(y)\, dS(y),
\end{equation}
an integral defined over a surface $\Gamma$,
where $\mathcal{K}$ is a free-space Green's function
for the underlying PDE, and
$\mu : \Gamma \to \mathbb{C}$ is a surface density function.
What makes this task
challenging is a combination of requirements for a
quadrature scheme, including the ability to handle singularities and near-singularities,
complex geometries $\Gamma$, on-surface and
near-surface evaluation, and support for simultaneous evaluation at a large
number of target points with low algorithmic complexity.

The potential that results from discretizing~\eqref{eqn:slp} with
a smooth high-order quadrature rule is a \emph{point potential} of the form
\begin{equation}
  \label{eqn:pt-pot}
  \ptpot(x_i) = \sum_{j=1}^N w_j \kernel(x_i, y_j)
  \quad
  (i=1,\dots,M),
\end{equation}
where $\{x_i\}_{i=1}^M \subseteq \mathbb{R}^3$ is a set of targets,
$\{y_i\}_{i=1}^N \subseteq \Gamma$ a set of sources, and $\{ w_i \}_{i=1}^N
\subseteq \mathbb{R}$ is a set of source weights. While $\ptpot$ is a generally
accurate approximation of $\mathcal{S}\mu$ for target points away from the source
geometry $\Gamma$, accuracy substantially decreases in the region close to
$\Gamma$ or on $\Gamma$ itself. A recent approach to resolving this problem is
based on Quadrature by Expansion (QBX,~\cite{klockner:2013:qbx}). By using the
broadly applicable assumption that the underlying kernel is a (locally) analytic
function of the target point $x$, and leveraging the strengths of smooth
high-order quadrature rules for one- and two-dimensional functions, QBX achieves
high-order accuracy without sacrificing generality. The key idea behind
QBX is that, to extend the applicability of smooth high-order quadrature to all
points $x \in \mathbb{R}^3$, the potential~\eqref{eqn:pt-pot} can be expanded in
a local expansion about a center $c \in \mathbb{R}^3 \setminus \Gamma$, and this
expansion may recover the value of the potential by analytic continuation in
regions where smooth quadrature is not adequate.  A potential of the
form~\eqref{eqn:pt-pot} can be evaluated in $O(N + M)$ time by way of the Fast
Multipole Method (FMM,~\cite{carrier:1988:adaptive-fmm}). Furthermore, as a
by-product of evaluation, the FMM forms local
expansions of the point potential covering the entire computational
domain. Since QBX expansions are precisely local expansions
of~\eqref{eqn:pt-pot} at appropriately chosen centers, the FMM therefore
provides a suitable path towards an `acceleration' strategy for reducing the
algorithmic complexity of QBX, by translating a suitably-chosen local expansion
to the QBX center. Care must be taken to ensure accurate combination of this
far-field potential with a suitable QBX-mediated near-field contribution.

The recent contribution~\cite{gigaqbx3d} combines a version of QBX, termed `global' QBX,
with the FMM for the evaluation of layer potentials in three
dimensions. This scheme, termed `Geometric Global Accelerated QBX', or `\algbrand',
in~\cite{gigaqbx2d,gigaqbx3d}, carefully controls the error introduced by
the FMM acceleration by enforcing strict geometric separation criteria
between intermediate multipole or local expansions
and QBX local expansions. As a consequence of these separation criteria,
the size of the near-field of a QBX center increases
if compared with a scenario where the QBX center is treated as a particle in a `point'-based FMM\@.
Empirically, the dominant cost of the scheme typically appears to come from
computing the spherical harmonic coefficients of QBX expansions from source
points in the near-field. A second important cost is the conversion of near-field
multipole expansions to QBX expansions. The primary purpose of this paper
is to reduce these two costs.

\begin{figure}
  \centering
  \begin{tikzpicture}[z={(0.1,-0.4)},scale=2.5]
    \draw[dashed] (-1, 0, -1/\sqrtthree)
      -- (1, 0, -1/\sqrtthree)
      -- (0, 0, 2/\sqrtthree)
      --  (-1, 0, -1/\sqrtthree);

    \foreach \srcX/\srcY in {
      0.000/0.000,
      0.787/-0.454,
      0.000/0.908,
      0.301/-0.454,
      -0.787/-0.454,
      0.243/0.488,
      0.544/-0.034,
      -0.243/0.488,
      -0.301/-0.454,
      -0.544/-0.034} {
      \draw[->,color=gray!40]
        (\srcX, 0, \srcY) node[fill=qbxcolor,coord] {}
        -- (\srcX, 0.5, \srcY) node[fill=qbxcolor,coord] {};
    }
    \DrawSphere[qbxcolor]{0}{0.5}{0}{0.5};

    \coordinate[fill=none,coord] (s1) at (0, 0, 0);
    \coordinate[fill=none,coord] (qbx) at (0, 0.5, 0);

    \draw (-1, 0, 0.7) node{stage-1 node} edge[->,out=0,in=180] (s1);
    \draw (-1, 1.1, 0.7) node{QBX center} edge[->,out=0,in=180] (qbx);
  \end{tikzpicture}
  \caption{Depiction of QBX center placement for a discretization over a
    triangular element. The QBX centers for on-surface evaluation are spawned in
    the normal direction at the discretization nodes
    termed `stage-1' nodes in~\cite{gigaqbx3d}, which serve as targets
    for on-surface evaluation.
    See also Section~\ref{sec:mesh-processing}.}%
  \label{fig:stage1}
\end{figure}
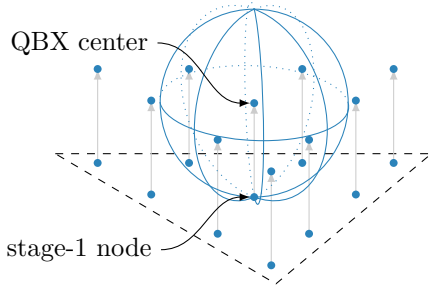

Siegel and Tornberg~\cite{siegel2017local} recently proposed using
\emph{target-specific} expansions to reduce the cost of forming QBX
expansions. Target-specific expansions are based on leveraging information about
both the source and the target to reduce the number of terms in the
expansion. Specifically, for expansions in spherical harmonics, they are based
on rewriting the expression for a local expansion using the
addition theorem for Legendre polynomials, an analytical tool previously applied to QBX in the context of
quadrature estimates~\cite{afklinteberg:2016:quadrature-est}. The number of
terms in the expansion is reduced in the polynomial order from $(p +
1)^2$ to $p + 1$, with a corresponding reduction in computational
effort, where $p$ is the expansion order.
Furthermore, the use of these expansions incurs no additional error
because they are based on a mathematical identity. The expansions are termed
`target-specific' because they do not separate the influence of the source and
target in the way that spherical harmonic expansions do, which means
that a different set of expansion `coefficients' is needed for each
source-target pair. Because of the need to recompute the coefficients for each
target, target-specific expansions cannot generally be used as a replacement for
spherical harmonic expansions in an FMM\@. (An exception to this is Anderson's
FMM~\cite{anderson:1992:fmm}, in which the expansion of potentials
is based on the Poisson integral formula, which is evaluated as a
target-specific expansion.)

Siegel and Tornberg introduce the numerical use of target-specific expansions in
the context of a scheme termed `local' QBX\@. The primary difference between
`global' QBX (such as the scheme in this paper) and `local' QBX is that, in
global QBX, the entire potential is conveyed to a target by a QBX expansion (hence `global' QBX),
whereas in local QBX only the potential due to sources in a
neighborhood of the target comes through the QBX
expansion, which is more akin to conventional FMM-accelerated
quadrature schemes employing local correction.
Compared with global QBX, local QBX features relatively
straightforward integration with the FMM due to relying only on `point'
evaluations for the accelerated part. A second advantage of local QBX is
that the placement of expansion centers is subject to less stringent geometric
requirements, leading to higher efficiency in areas where the geometry is close
to touching or highly irregular. Despite these advantages, for
achieving the same level of accuracy as global QBX, current versions of local QBX appear to require
higher quadrature order
and oversampling in the local neighborhood of a target. This appears to be due to a form of
truncation error not present in global QBX---specifically, error introduced
when matching the transition between the QBX-mediated local neighborhood and the
`point' field from far-away sources.

The main observation in this paper is that target-specific expansions,
originally developed for local QBX and used a post-processing step in
combination with a `point' FMM, can also be used within the context of global
QBX, and, more specifically, within the \algbrand\ FMM\@.  An
approximation to the QBX expansion $\Lqbx{c}(t)$ in the \algbrand~FMM at a
target $t \in \mathbb{R}^3$ associated with an expansion center $c \in
\mathbb{R}^3$ is formed from three parts (see Section~\ref{sec:formal-statement} for notation)
\begin{equation}
  \Lqbx{c}(t) = \Lqbxfar{c}(t) + \LqbxW{c}(t) + \Lqbxnear{c}(t),
  \label{eq:potential-pieces}
\end{equation}
where the portion $\Lqbxnear{c}$ is
mediated by direct formation of the expansion due to the sources in the
near-field of the box `owning' the center $c$, and the quantities $\Lqbxfar{c}$
and $\LqbxW{c}$ are obtained through the various expansion translations
present in the Fast Multipole Method. Unlike local QBX,
analytical truncation error estimates for global QBX only apply to the
combined expansion $\Lqbx{c}$ and not the sub-expansions. Thus, while no
truncation bounds are known for the individual terms in~\eqref{eq:potential-pieces},
truncation bounds do hold on the overall sum, and the partitioning
of the potential into sub-expansions has no effect on truncation error. In the
context of the present work (and in keeping with the cost argument above)
we use target-specific expansions to evaluate the component $\Lqbxnear{c}$.

The observation that a component of the FMM field can be mediated with
target-specific expansions applies conceptually to not only the
\algbrand~variant of global QBX but to others as well such as the ones
in~\cite{rachh:2017:qbx-fmm,rahimian:2017:qbkix}, or three-dimensional versions thereof.
It is however likely that
target-specific expansions will have greater cost impact on \algbrand~due
to the larger proportion of direct interactions.

An especially important use case, and one for which target-specific expansions
can excel in cost compared with (`target-independent') spherical harmonic
expansions, is \emph{on-surface evaluation}. In our current treatment of
on-surface evaluation, each on-surface target uses a different QBX center
(Figure~\ref{fig:stage1}), meaning
that no QBX expansion is typically evaluated at more than one target.
As such, there is no advantage to be had from the source-target separation
of variables in spherical harmonic expansions, and when there is only one target per center, forming
and evaluating spherical harmonic expansions is more expensive than target-specific expansions.

A second observation we make in this paper is that, while replacing near-field
evaluations with target-specific expansions results in cost improvements,
adjusting parameters to the algorithm expectedly leads to further opportunities
for cost reduction. By design, the FMM is tasked with making choices between
whether to evaluate an interaction directly or mediate it via expansions. Assuming one is presented with a set of evaluation strategies which meet required accuracy tolerances,
a standard method for minimizing computational cost is to use
thresholds based on particle counts to decide which strategy
to use. Because target-specific expansions make direct evaluations less expensive,
adjusting the thresholds in such a way as to shift a larger fraction of the work onto
direct interactions reduces the overall cost of the FMM\@. To aid in this
\emph{rebalancing} of the FMM, in this paper we develop a cost model that
models the number of floating point operations in the \algbrand~FMM\@. With
appropriate fitting of empirical per-stage calibration factors, this model is
able to approximate the total computational time used by the algorithm with high
accuracy. This in turn provides a reproducible cost measure, used here for
wall-time independent balancing and the reporting of cost and scalability
results.

Much analytical
modeling work for optimizing for FMM cost is restricted to the case of uniform
distributions (e.g.~\cite{cheng:1999:3d-fmm,gumerov-helmholtz-fmm,singer:1995:fmm,zhao:1991:fmm}),
which makes
it inapplicable to layer potential evaluation, since the particle distributions
arising from surface discretizations for layer potentials are not uniform in the volume.  Our
work differs from the approaches for uniform distributions by making use of more
information from the geometry in order to give a precise prediction of cost. A
similar approach to FMM cost modeling based on approximating the number of
floating point operations through direct inspection of the geometry is used
in~\cite{anderson:1992:fmm,hu:1996:accuracy}, though the details of the modeled algorithms differ
substantially from~\algbrand. Other work on cost models for
nonuniform particle distributions includes the
contributions~\cite{pouransari:2015:adaptive-fractal-fmm,agullo:hal-01474556,nabors:1994:fmm}.
In~\cite{pouransari:2015:adaptive-fractal-fmm}, a model is developed to optimize
for FMM parameters when the particle distribution is a fractal set.
In~\cite{agullo:hal-01474556}, an empirical model is developed for
predicting the cost of a particular task-based FMM implementation on arbitrary
particle distributions. The contribution~\cite{nabors:1994:fmm} discusses
general conditions on the distribution under which linear scaling may be expected.
In addition to optimization for cost, models have been applied to aid in solving
the problem of distributing work among heterogeneous
systems~\cite{Choi:2014:CGH:2588768.2576787} or predicting execution
characteristics taking into account both computation and memory
accesses~\cite{Chandramowlishwaran:2012:BAT:2312005.2312039}.

In summary, we present the following contributions in this paper.

\begin{itemize}
  \item We describe how to use target-specific expansions (TSQBX) inside the
    global QBX (\algbrand) FMM, a technique previously only used in `local' QBX,
    to reduce the cost of the near-neighborhood interactions.
  \item We present a cost model for the \algbrand~algorithm on a shared
    memory system, which predicts the running time of the \algbrand~FMM with
    very high accuracy.
  \item We demonstrate a 1.7--3.3$\times$ reduction in modeled cost using TSQBX on
    test cases with complex unstructured geometries
    for the Laplace and Helmholtz equations.
\end{itemize}

The organization of this paper is as follows. Section~\ref{sec:background}
describes background material pertaining to QBX and the \algbrand~FMM\@. In
Section~\ref{sec:algorithm}, we present for completeness a full statement of
the~\algbrand\ algorithm using target-specific expansions. In
Section~\ref{sec:results}, we present a study pertaining to the cost
impact of target-specific expansions within the \algbrand~FMM, making concluding
remarks in Section~\ref{sec:conclusion}. Appendix~\ref{sec:ts-expansions}
presents a derivation of target-specific expansions for various kernels.
Appendix~\ref{sec:software} describes how to obtain the software used in this paper.

\section{Background}%
\label{sec:background}

\subsection{Layer Potentials}

For the sake of exposition, we consider the solution of the exterior Neumann
problem in three dimensions, for a smooth, bounded, simply or multiply-connected
domain $\Omega$ with boundary $\Gamma$. For continuous Neumann data~$g$, the
boundary value problem is to find a function~$u: \mathbb{R}^3 \setminus \Omega
\to \mathbb{R}$ such that
\begin{DIFnomarkup}
\begin{alignat}{2}
  \triangle             u(x) &= 0 & \quad & x \in \mathbb{R}^3 \setminus \Omega, \nonumber \\
  \lim_{h\to 0^{+}}
  \hat \nu(x) \cdot \nabla u(x + h \hat \nu(x))
                          &= g(x) & \quad & x \in \Gamma, \nonumber \\
  \lim_{|x| \to \infty} u(x) &= 0.
\end{alignat}
\end{DIFnomarkup}
The notation~$\hat \nu(x)$ refers to the outward-facing unit normal vector at $x$.
The method under consideration here lends itself to the solution of a considerably broader family of boundary value
problems.

We represent the solution of the problem by means of layer potentials.  In the
remainder of this paper, we will use $\norm{ \, \cdot \,}$ to denote the
Euclidean ($\ell^2$) norm unless otherwise specified.  Recalling the Green's function for
the Laplace equation,
\begin{equation}
  \label{eqn:laplace-green-function}
  \mathcal{G}(x, y) = (4 \pi)^{-1} \norm{x - y}^{-1},
\end{equation}
we consider in this section the single-layer potential~\eqref{eqn:slp}
$\mathcal{S}\mu$ with kernel $\kernel = \mathcal{G}$ and we introduce the following
layer potential $\mathcal{S}'\sigma$ with density function~$\sigma: \Gamma \to
\mathbb{R}$:
\[
  (\mathcal{S}' \sigma)(x) \coloneqq \int_\Gamma
  \frac{\partial \mathcal{G}(x, y)}{\partial \hat{\nu}(x)} \sigma(y) \, dS(y).
\]
With the aid of these operators, we represent the solution~$u$ as
\[
  u \coloneqq \mathcal S \mu
\]
using an unknown density $\mu : \Gamma \to \mathbb{R}$. By differentiation
under the integral sign, any function of the form~$\mathcal{S} \mu$
satisfies the Laplace PDE in the exterior of~$\Omega$, and by considering the
asymptotic behavior of the Green's function~$\mathcal{G}$, it also automatically
satisfies the far-field decay conditions.

We enforce the Neumann boundary conditions on this representation as
follows. The classical \emph{jump
  relations}~\cite[Thm.~6.19]{kress:2014:integral-equations} imply that the normal
derivative of $\mathcal{S} \mu$ is discontinuous across the boundary, in the
sense that
\[
  \lim_{h\to 0^{\pm}}
  \hat \nu(x) \cdot \nabla u(x + h \hat \nu(x)) = \mathcal{S}' \mu(x) \mp
  \frac{1}{2} \mu(x)\quad(x\in\Gamma).
\]
Thus, to satisfy the boundary condition, $\mu$ needs to solve the following
second kind integral equation:
\begin{equation}
  \label{eqn:exterior-neumann-ie}
  g = \left( \mathcal{S}' - \frac{1}{2} \right) \mu.
\end{equation}
With the help of the Fredholm theory for second kind integral equations one obtains that the solution to
this equation exists, is unique, and is continuously dependent on
$g$~\cite[Thm.~6.28,~6.30]{kress:2014:integral-equations}. Under a suitable
discretization, this equation provides the basis for numerical methods for
solution of the exterior Neumann problem.

\subsection{High-Order Quadrature for Smooth Functions}%
\label{sec:quadrature}

The primary concern of this paper is the numerical evaluation of layer
potentials such as the single layer potential~\eqref{eqn:slp} anywhere in
$\mathbb{R}^3$, including near or on the surface $\Gamma$. A natural but
ultimately deficient approach to this problem uses smooth composite
quadrature. In this approach, the surface is assumed to be tessellated into $K$
disjoint surface elements
\[ \Gamma = \bigcup_{k=1}^K \Gamma_k. \]
Each element $\Gamma_k$ is parametrized by a smooth mapping function
$\mapping_k: E \to \mathbb{R}^3$ from a reference element $E$ in the
plane. Then, with the use of an $M$-point quadrature rule defined over the
reference element with weights $\{w_i\}_{i=1}^M$ and nodes $\{y_i\}_{i=1}^M$,
the single layer potential admits the approximation
\begin{align*}
  (\mathcal{S} \mu)(x)
  &= \sum_{k=1}^K
  \iint_E \mu(\mapping_k(y)) \cdot \mathcal{G}(x, \mapping_k(y))
  \cdot \left| (\partial_{e_1} \mapping_k \times
  \partial_{e_2} \mapping_k)(y) \right|
  \, dS(y) \\
  &\approx
  \sum_{k=1}^K \sum_{i=1}^M w_i \cdot \mu(\mapping_k(y_i))
  \cdot \mathcal{G}(x, \mapping_k(y_i)) \cdot \left| (\partial_{e_1} \mapping_k \times
\partial_{e_2} \mapping_k)(y_i) \right|.
\end{align*}
For points $x$ far from the surface $\Gamma$, the integrand is a smooth function
favorable to numerical treatment with conventional quadrature rules for smooth
functions on plane regions. However, for $x$ near to the surface, the quadrature
error is at worst unbounded due to near-singularity of the Green's function
under the integrand, which necessitates a massive increase in quadrature order
to resolve. This effectively prevents the practical use of traditional high-order
composite quadrature for these integrands as an evaluation strategy in a neighborhood of the surface,
where the size of the neighborhood depends on the size of the element, the
quadrature order, and the desired accuracy~\cite{klockner:2013:qbx}.

The key insight in QBX as a quadrature scheme is that $\mathcal{S} \mu$ is
an analytic function on $\mathbb{R}^3 \setminus \Gamma$. This fact can be
exploited to recover quadrature accuracy for $x$ near the surface via analytic
continuation, at a lower cost than smooth quadrature.



\subsection{Expansion of Potentials in Spherical Harmonics}%
\label{sec:sph-harm}

We recall the addition theorem for the Laplace potential in three
dimensions. Let $a, b \in \mathbb{R}^3$ with $0 <
\norm{a} < \norm{b}$, and let $\gamma_{ab}$ be the angle (about the origin) between $a$ and $b$. By
expanding $\mathcal{G}$ in a binomial series~\cite{kellogg}, the
potential can be represented as
\begin{equation}
  \label{eqn:laplace-potential-expansion}
  \mathcal{G}(a, b) =
  \frac{1}{4 \pi}
  \sum_{n=0}^\infty \frac{\norm{a}^n}{\norm{b}^{n+1}} P_n(\cos \gamma_{ab}).
\end{equation}
The function $P_n$ is the Legendre polynomial of degree $n$. The term $P_n(\cos
\gamma_{ab})$ may be further expanded in a series of spherical harmonics.
Let $a$ and $b$ be written in polar and azimuthal spherical coordinates
$(\theta_a, \phi_a)$ and $(\theta_b, \phi_b)$ respectively, i.e.\ with $\theta =
\cos^{-1}(z/r), \phi = \operatorname{atan2}(y, x).$ The identity known as the
spherical harmonic addition theorem, or the addition theorem for Legendre polynomials,
says that~\cite{siegel2017local}
\begin{equation}
  \label{eqn:legendre-addition-identity}
  P_n(\cos \gamma_{ab}) =
  \frac{4 \pi}{2n + 1}
  \sum_{m=-n}^n Y_n^m(\theta_a, \phi_a) Y_n^{-m}(\theta_b, \phi_b).
\end{equation}
The spherical harmonics $Y_n^{m}$, $m, n \in \mathbb{N}_0$, $|m| \leq n$, are defined
(following~\cite{siegel2017local}) as
\begin{equation}
  \label{eqn:sph-harm}
  Y^m_n(\theta, \phi) \coloneqq \sqrt{\frac{2n+1}{4 \pi}
    \frac{\left(n-|m|\right)!}{\left(n+|m|\right)!}}
    \cdot P^{|m|}_n(\cos \theta) e^{i m \phi},
\end{equation}
where $P^m_n$ is the associated Legendre function of order $m$ and degree $n$.
Substituting~\eqref{eqn:legendre-addition-identity}
into~\eqref{eqn:laplace-potential-expansion}, we obtain that for $0 \leq \norm{a} <
\norm{b}$,
\begin{equation}
  \label{eqn:laplace-identity}
  \mathcal{G}(a,b) = \sum_{n=0}^\infty \frac{1}{2n + 1}
  \frac{\norm{a}^n}{\norm{b}^{n+1}} \sum_{m=-n}^n Y_n^m(\theta_a, \phi_a)
  Y_n^{-m}(\theta_b, \phi_b).
\end{equation}
This series allows us to expand the Green's function as follows. For a given
choice of source point $b$ and expansion center $c$, one can define local
expansion coefficients as a doubly-indexed sequence
\begin{equation}
  \label{eqn:local-coeffs}
  L_n^m \coloneqq \frac{1}{2n+1}
  \frac{Y_n^{-m}(\theta_{b-c}, \phi_{b-c})}{|b-c|^{n+1}}
\end{equation}
for integer $|m| \leq n$ and $n \in \mathbb{N}_0$. Then the Green's function
evaluated at a target $a$, with $|a - c| < |b - c|$, can be written as
\begin{equation}
  \label{eqn:local-expansion}
  \mathcal{G}(a, b) =
  \sum_{n=0}^\infty \sum_{m=-n}^n L_n^m |a-c|^n Y_n^m(\theta_{a-c}, \phi_{a-c}).
\end{equation}
A $p$-th order expansion is one in which the series~\eqref{eqn:local-expansion}
is truncated to the first $p + 1$ terms.

\subsection{Quadrature by Expansion}

Using the local expansion of the potential as an analytical tool, we are ready
to discuss the fundamentals of QBX\@. Among the current work on QBX, two major
variants of the scheme have been considered by various authors: `global' and
`local.' The variant that is the subject of this paper is a global scheme
that borrows some ideas from recent work on local QBX. Hence, we review both of these
variants in this section. QBX may be understood as a discretization involving
two inter-related stages: formation of a truncated local expansion and smooth
quadrature.

\paragraph*{First Stage: Formation of a Truncated Local Expansion.}
With reference to the source surface $\Gamma$, for each target point $t$
close to or on the surface, this stage chooses a point $c$ to act as an expansion
center. The distance $\lvert t - c \rvert$ from $c$ to $\Gamma$ is called the \emph{expansion
  radius}. Using the selected center, this stage forms a local expansion about
$c$ to mediate the potential $\mathcal{S}\mu$.

In \emph{global} QBX, the potential $\mathcal{S}\mu$ due to the entire source geometry is expanded about
$c$. A sufficient condition for the convergence of the expansion is that $\lvert t - c \rvert \leq
\mathrm{dist}(c, \Gamma)$. By applying~\eqref{eqn:local-coeffs}, one defines QBX
coefficients via the integrals
\begin{equation}
  \label{eqn:qbx-coeff}
  (L_\mathrm{global})_n^m(c) \coloneqq
  \frac{1}{2n + 1} \int_\Gamma \frac{Y_n^{-m}(\theta_{y-c},
  \phi_{y-c})}{\norm{y-c}^{n+1}} \mu(y) \, dS(y).
\end{equation}
Then, fixing an expansion order $p \in \mathbb{N}_0$, the coefficients
$(L_\mathrm{global})_n^m(c)$, for $|m| \leq n$ for all $n \leq p$, are used to
approximate the single layer potential in a series centered at $c$
\begin{equation}
  \label{eqn:qbx-series}
  \qbxsinglelayercont{p} \mu(t) \coloneqq
  \sum_{n=0}^p \sum_{m=-n}^{n} (L_\mathrm{global})_n^m(c) \norm{t-c}^n
  Y_n^m(\theta_{t-c}, \phi_{t-c}).
\end{equation}

In contrast, local QBX mediates only a part of $\mathcal{S}\mu$
due to source geometry in  a neighborhood of the target via a local expansion.
In \emph{local QBX},
formation of the expansion starts with the splitting
\[
  (\mathcal{S}\mu)(t) =
    \left(\mathcal{S} \left. \mu \right|_{\Gamma_{\mathrm{local},t}}\right)(t)
    + \left(\mathcal{S} \left. \mu \right|_{\Gamma \setminus \Gamma_{\mathrm{local},t}}\right)(t),
\]
where the region $\Gamma_{\mathrm{local},t} \subseteq \Gamma$ is chosen in such a way as
to include the nearly singular/singular portion of the integrand. We therefore have
\begin{align*}
  \left(\mathcal{S} \left.\mu\right|_{\Gamma_{\mathrm{local},t}}\right)(t)
    &= \int_{\Gamma_{\mathrm{local},t}} \mathcal{G}(t, s) \mu(s) \, dS(s), \\
  \left(\mathcal{S} \left.\mu\right|_{\Gamma \setminus \Gamma_{\mathrm{local},t}}\right)(t)
    &= \int_{\Gamma \setminus \Gamma_{\mathrm{local},t}} \mathcal{G}(t, s) \mu(s) \, dS(s).
\end{align*}
One then should ensure that
$\lvert t - c \rvert \leq \mathrm{dist}(c, \Gamma_{\mathrm{local},t})$.
The QBX coefficients are defined via
integrals
\begin{equation}
  \label{eqn:local-qbx-coeff}
  (L_{\mathrm{local},t})_n^m(c) \coloneqq
  \frac{1}{2n + 1} \int_{\Gamma_{\mathrm{local},t}}
  \frac{Y_n^{-m}(\theta_{y-c},
  \phi_{y-c})}{\norm{y-c}^{n+1}} \mu(y) \, dS(y),
\end{equation}
and the approximation to the potential is given by
\begin{equation}
  \label{eqn:local-qbx-series}
  \localqbxsinglelayercont{p}{t} \mu(t)
  \coloneqq   \sum_{n=0}^p \sum_{m=-n}^{n} (L_{\mathrm{local},t})_n^m(c) \norm{t-c}^n
  Y_n^m(\theta_{t-c}, \phi_{t-c})
  + \left(\mathcal{S} \left. \mu \right|_{\Gamma \setminus \Gamma_{\mathrm{local},t}}\right)(t).
\end{equation}

The difference between the expansions~\eqref{eqn:qbx-series}
or~\eqref{eqn:local-qbx-series} and the value of the layer potential is termed
the truncation error. The following result gives an accuracy estimate for the
expansion~\eqref{eqn:qbx-series} used by global QBX, for the case that $t$ is an
on-surface target. Bounds are available for the off-surface case through the
same analysis.
\begin{lemma}[QBX truncation error, based on~{\cite[Thm~3.1]{epstein:2013:qbx-error-est}}]%
  \label{lem:qbx-truncation-3d}%
  Suppose that $\Gamma$ is smooth, non-self-intersecting and let $r > 0$.
  Suppose that $\{x : \norm{x - c} \leq r \} \cap \Gamma = \{ t \}$.
  Then for each $p > 0$ and $\delta > 0$, a constant $M_{p,\delta}$
  exists such that
  \begin{equation}
    \label{eqn:truncation-estimate}
    \left|
    (\mathcal{S}\mu)(t) -
    \qbxsinglelayercont{p} \mu(t)
    \right| \\
    \leq M_{p, \delta} r^{p+1}
    \| \mu \|_{W^{3 + p + \delta, 2}(\Gamma)}.
  \end{equation}
\end{lemma}

Error estimates for local QBX evaluation~\eqref{eqn:local-qbx-series} can be
found in~\cite{siegel2017local}. The primary difference between the truncation
error in global QBX and local QBX is that for local QBX, there is a dependence
on the ratio $r/R$, where $r$ is the expansion radius and $R$ is the distance
between the expansion center $c$ and the nearest boundary point of the surface region
$\Gamma_{\mathrm{local},t}$. This error can be interpreted as the error
associated with the non-smooth transition between the expansion-mediated
contribution from the local neighborhood $\Gamma_{\mathrm{local},t}$ and the
non-expansion-mediated `far-field.'

\paragraph*{Second Stage: Smooth High-Order Quadrature.}

In the second stage, the QBX coefficients $(L_\mathrm{global})_n^m(c)$ or
$(L_{\mathrm{local},t})_n^m(c)$, and, in the local case, the `far-field'
integral are approximated with a quadrature rule, commonly
with smooth high-order `panel'-based quadrature discussed in
Section~\ref{sec:quadrature}. When the source points are close to the expansion center,
accuracy of the quadrature approximation depends strongly on the expansion radius.
Specifically, to maintain high-order
quadrature accuracy, the expansion center must be placed sufficiently far from
the surface element spawning the center as well as nearby elements,
where the critical distance depends on
the `size' of the element. As an example providing quantitative detail, the
following result, due to af~Klinteberg and
Tornberg~\cite{afklinteberg:2016:quadrature-est}, gives the asymptotic error for
the case of a smooth tensor product rule over a flat $2h \times 2h$ panel.
\begin{lemma}[QBX quadrature error,
    flat $2h \times 2h$ panel~{\cite[eqn.~(157)]{afklinteberg:2016:quadrature-est}}]%
  \label{lem:qbx-quadrature-3d}%
  Let $\Gamma = [-h,h] \times [-h,h] \times \{0\}$, and let $t = (x, y, 0) \in
  \Gamma$ be a target point and $c = (x, y, r)$ be the corresponding expansion
  center. Let $\mu: \Gamma \to \mathbb{R}$ be a density function defined on $\Gamma$.
  Assume that $\mu$ is smooth (see~\cite[Sec.~3.3]{afklinteberg:2016:quadrature-est} for more discussion concerning smoothness of $\mu$).
  Suppose that the coefficients $(L_\mathrm{global})_n^m(c)$ to the
  series~\eqref{eqn:qbx-series} up to order $p$
  are computed using a $q$-point Gauss-Legendre
  tensor product rule.
  Then a constant $C > 0$ exists, independent of
  $h$, $\mu$, $x$, $y$, $r$, $p$, and $q$, such that the
  quadrature-based approximation $\qbxsinglelayer{p} \mu(t)$ to the QBX expansion
  $\qbxsinglelayercont{p} \mu(t)$
  satisfies the following error bound
  asymptotically as $q \to \infty$:
  \[
     \left|
     \qbxsinglelayercont{p} \mu(t) -
     \qbxsinglelayer{p} \mu(t)
     \right|
     \lesssim
     C \frac{h}{q} \sum_{l=1}^p \frac{1}{\sqrt{l}} \left( \frac{4qre}{hl} \right)^l
     e^{-4qr/h} \| \mu \|_\infty.
  \]
\end{lemma}

When combined,
Lemmas~\ref{lem:qbx-truncation-3d}~and~\ref{lem:qbx-quadrature-3d} suggest that,
with careful control over the expansion radius, quadrature order, and element
size, QBX may be used to obtain high-order quadrature accuracy for on-surface
targets $t$ (and, completely analogously, for off-surface targets as well). The
two stages of QBX described in this section are interrelated, in that the choice
of quadrature discretization influences the choice of expansion radius. For
local QBX, an additional degree of freedom is the size of the local neighborhood
$\Gamma_{\mathrm{local},t}$, which has non-trivial implications on truncation
error and quadrature error.

The remainder of this paper focuses on the global scheme. The next section
discusses how to achieve the type of error control described in the previous
paragraph for global QBX, on user-supplied meshes of smooth geometries.

\subsection{Surface Discretization}%
\label{sec:mesh-processing}

\begin{figure}
  \centering
  \begin{tikzpicture}[z={(0.1,-0.4)},scale=2.5]
    \draw[dashed] (-1, 0, -1/\sqrtthree)
      -- (1, 0, -1/\sqrtthree)
      -- (0, 0, 2/\sqrtthree)
      --  (-1, 0, -1/\sqrtthree);

    \foreach \srcX/\srcY in {
      0.000/0.000,
      0.787/-0.454,
      0.000/0.908,
      0.301/-0.454,
      -0.787/-0.454,
      0.243/0.488,
      0.544/-0.034,
      -0.243/0.488,
      -0.301/-0.454,
      -0.544/-0.034} {
      \draw (\srcX, 0, \srcY) node[fill=qbxcolor, coord] {};
    }

    \foreach \srcX/\srcY in {
      0.412/-0.238,
      -0.412/-0.238,
      0.213/0.123,
      0.000/-0.246,
      0.000/0.476,
      -0.213/0.123,
      0.209/-0.509,
      -0.209/-0.509,
      0.545/0.074,
      -0.545/0.074,
      0.337/0.435,
      -0.337/0.435,
      0.617/-0.506,
      0.746/-0.281,
      -0.617/-0.506,
      -0.746/-0.281,
      0.130/0.787,
      -0.130/0.787,
      0.906/-0.523,
      -0.906/-0.523,
      0.000/1.046} {
      \draw[color=gray!40] (\srcX, 0, \srcY) node[fill=gray!50, coord] {}
        .. controls (0, 0.3, 0) .. (0,0.45,0);
    }

    \draw[->,color=gray!50] (0,0.45,0) -- (0,0.5,0)
      node[fill=qbxcolor, coord] {};

    \DrawSphere[qbxcolor]{0}{0.5}{0}{0.5};

    \coordinate[coord] (s1) at (0,0,0);
    \coordinate[coord] (s2q) at (0,0,0.476);
    \coordinate[coord] (qbx) at (0,0.5,0);

    \draw (-0.8, 0, 1) node {stage-2 quad.\ node} edge[->,out=0,in=270] (s2q);
    \draw (-0.9, 0, 0.5) node{stage-1/2 node} edge[->,out=0,in=180] (s1);
    \draw (-1, 1, 0.7) node{QBX center} edge[->,out=0,in=180] (qbx);
  \end{tikzpicture}
  \caption{Depiction of a stage-1 discretization, stage-2 discretization,
    and a stage-2 quadrature discretization on a triangular element.
    In this example the stage-1 and stage-2 discretizations coincide; in general
    they may be different.
    The stage-2 quadrature discretization provides the quadrature
    nodes for the local expansion of the potential formed at the QBX
    centers. See also Figure~\ref{fig:stage1}.}%
  \label{fig:stage2-quad}
\end{figure}
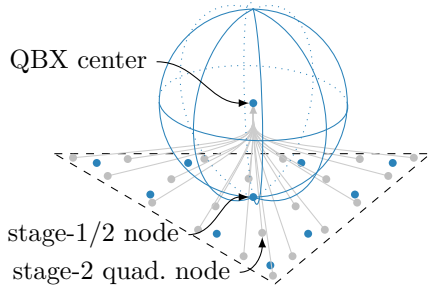

In the remainder of this paper, we will assume that the surface discretization
is an unstructured triangular surface mesh given as the union of images of a
triangular reference element under a polynomial mapping. Our specific choice of
basis for reference nodes/degrees of freedom
follows~\cite{vioreanu_spectra_2014}, and for quadrature nodes is taken
from~\cite{xiao_numerical_2010}. A density on the discretization is represented
by its values at the mapped counterparts of quadrature or interpolation nodes
given on a reference element.

As part of the \algbrand~algorithm in three dimensions,
the work in~\cite{gigaqbx3d} describes an algorithm for preparing an arbitrary
smooth geometry to serve as a surface source/target discretization for
applying QBX\@. We briefly review its main steps for the
benefit of the reader. More details, including considerations for efficient
implementation, may be found there.

The mesh processing algorithm for \algbrand~involves four versions of the
surface discretization, connected via interpolation operators that bring density
values from one version to a subsequent refined version. First, the unmodified mesh is
received from an external mesh generator. Second, the mesh is refined by
iterative bisection of the elements~\triforce{} to avoid conflicts of the QBX
expansion balls with other source geometry, to obtain the version called the
`stage-1 discretization'. Expansion centers are placed in the normal directions
at the stage-1 interpolation nodes (which become on-surface targets) at a
distance proportional to the element size (Figure~\ref{fig:stage1}), and remain
fixed in the subsequent discretizations. Third, the mesh is refined with
iterative bisection to assure sufficient quadrature resolution for all
interactions (called the `stage-2 discretization').  Lastly, the quadrature
nodes in the mesh are oversampled (by increasing the quadrature order) to ensure accurate evaluation of the QBX
coefficient integrals~\eqref{eqn:qbx-coeff}, producing the `stage-2 quadrature
discretization' (Figure~\ref{fig:stage2-quad}), providing the source points for
the `point potential' approximation to the layer potential. The latter
discretization is suitable as input to an algorithm for the evaluation of point
potentials such as the FMM\@.

\subsection{FMM Acceleration}%
\label{sec:fmm-acceleration}

The focus of the contributions~\cite{gigaqbx2d, gigaqbx3d} is the accelerated
evaluation of the quadrature for layer potentials. This is done using an
appropriately modified version of the Fast Multipole Method (FMM).
Algebraically, QBX may be regarded as the evaluation of the \emph{local
  expansion} of the point potential~\eqref{eqn:pt-pot} due to the quadrature
nodes. Because such local expansions may be viewed as the output of the FMM, a
natural approach to modifying the FMM for QBX is to use the FMM to form local
expansions of the potential at the QBX expansion centers. The first practical
implementation of this approach, as described in~\cite{rachh:2017:qbx-fmm}, is
not backed by error estimates and does not achieve the same level of accuracy for
a given FMM expansion order as the `point' FMM, though an empirically
determined increase in FMM order can recover accuracy, at some expense. As explored
in detail in~\cite{gigaqbx2d}, the geometrical root cause of this loss of
accuracy is that the FMM separation constraints for accurate evaluation of
`point' potentials are not strong enough to prevent inaccurate contributions
from entering the QBX expansion.

The main modification to the QBX FMM in~\cite{gigaqbx2d, gigaqbx3d} permits
targets, such as QBX expansion balls, to be `sized'.
Similar to point targets, sized targets have a near-field that disallows certain nearby
sources from using expansion mediation, permitting for analytical accuracy bounds to be established.
This is accounted for by only allowing sized targets to
protrude beyond their containing boxes by at most a given factor relative to the
box size, called the \emph{target confinement factor}. If a sized target cannot
fit in a child box, it remains in the parent box. These changes require a
careful reworking of the definitions of a number of aspects of the classical
Fast Multipole algorithm for accuracy and
scalability of the resulting method. The resulting method was termed
\emph{GIGAQBX}---for `Geometric Global Accelerated QBX'---in~\cite{gigaqbx2d, gigaqbx3d}.

The redefinition of the near-field of a sized target entails a larger number of
`direct' FMM interactions at the QBX centers. Assuming the size of the QBX near-field
remains bounded, this does not threaten the
theoretical scaling of the algorithm. However, a practical implication of this
is that, in three dimensions, direct interactions may take a large portion
of time due to the high cost of expansion formation. The main contribution of this paper
consists of an approach for reducing this cost, discussed next.

\subsection{Target-Specific Expansions}%
\label{sec:ts}

\begin{figure}
  \centering
  \begin{tikzpicture}[z={(-0.09,-0.09)},scale=1.5]
    \DrawSphere[qbxcolor]{0}{0}{0}{1}

    \coordinate[fill,coord] (c) at (0,0,0);
    \coordinate[fill,coord] (s) at (-1.5,-0.75,1);
    \def\tc{1/sqrt 3}
    \coordinate[fill,coord] (t) at (-\tc,\tc,\tc);

    \draw[<->] (s) -- (c);
    \draw[<->, label] (c) -- (t);
    \draw (c) node[right] {$c$};
    \draw (s) node[left] {$s$};
    \draw (t) node[left] {$t$};


    \begin{scope}[canvas is plane={O(0,0,0)x(-\tc,\tc,\tc)y(-0.549,-0.798,0.249)}]
      \draw (0.4, 0) arc (0:58.83:0.4) node[midway, left] {$\gamma$};
    \end{scope}
  \end{tikzpicture}
  \caption{Components of a target-specific QBX expansion.}%
  \label{fig:tsqbx}
\end{figure}
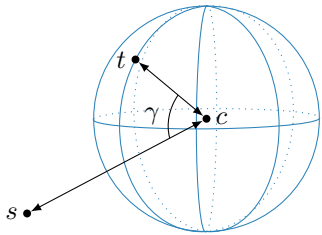

Siegel and Tornberg~\cite{siegel2017local} observe
that the series expansion of the Green's
function~\eqref{eqn:laplace-potential-expansion} provides a way to accelerate
the formation of QBX expansions in certain circumstances. Consider the computational
problem of evaluating the local expansion due to $n_s$ sources at $n_t$ targets. Assume
for simplicity that we use a single expansion center mediating the expansion
of the entire potential. We describe two approaches to this problem.

In the first approach, using the formula~\eqref{eqn:laplace-identity}, for each
source point $s$ we compute $(p + 1)^2$ intermediate local coefficients, which
takes $O((p+1)^2)$ time per source with well-known recurrences.  After combining
local coefficients additively to obtain $(p + 1)^2$ final coefficients, we
evaluate the local expansion at each target point, which costs $O((p + 1)^2)$ at
each target, again using well-known recurrences. It follows that this approach
requires $O((n_s + n_t)(p+1)^2)$ work.

In the second approach, we use
formula~\eqref{eqn:laplace-potential-expansion}. Recall that this formula
implies that the $p$-th order local expansion of the potential due to $s$ about
the center $c \in \mathbb{R}^3$, with coefficients $L_n^m$, satisfies
\begin{equation}
  \label{eqn:tsqbx}%
  \sum_{n=0}^p \sum_{m=-n}^n L_n^m |t-c|^{n} Y_n^m(\theta_{t-c}, \phi_{t-c})
  =
  \frac{1}{4\pi} \sum_{n=0}^p \frac{\norm{t-c}^n}{\norm{s-c}^{n+1}} P_n(\cos \gamma),
\end{equation}
where $\gamma$ is the angle between $s - c$ and $t - c$. A depiction of the
geometrical situation is given in Figure~\ref{fig:tsqbx}. The
quantity~\eqref{eqn:tsqbx} requires $p + 1$ summation terms and, using
recurrences for the Legendre polynomials $P_n$, can be evaluated in $O(p + 1)$
time. It must be evaluated once per source/target pair. It follows that the
total cost of the second approach is $O(n_s n_t (p + 1))$.

The first approach is the one used internally within the FMM\@. It has the advantage
of scaling linearly in the number of sources and targets. The disadvantage is
the high cost of $O((p+1)^2)$ operations per particle.

The second approach, while not scaling linearly in the number of particles,
is computationally advantageous over the first approach if the number of sources
or targets is small. Such a situation arises in on-surface evaluation for QBX
when there is one target per center (Figure~\ref{fig:stage1}), a common evaluation pattern.

The second approach is termed \emph{target-specific} because, unlike the first
approach, the `local coefficient' $P_n(\cos \gamma) / \norm{s - c}^{n+1}$
depends on the target through the angle $\gamma$.  Because of this, it is
generally unsuitable for use within a point FMM, which is premised on using
expansions that separate the influence of the source and the targets. However,
target-specific expansions are suitable for use within the \algbrand\ FMM
whenever a QBX local expansion is formed directly from sources. Furthermore,
analogous formulas for target-specific QBX are available for different kernels
derived from the Laplace and Helmholtz kernels, making this approach general. We
give some of these in Appendix~\ref{sec:ts-expansions}.

In the next section, we discuss the incorporation of target-specific QBX into
the \algbrand\ algorithm.

\section{Algorithm}%
\label{sec:algorithm}

The algorithm in this section is based on the \algbrand~algorithm in three
dimensions~\cite{gigaqbx3d}. For the benefit of the reader familiar with the
version in~\cite{gigaqbx3d}, we briefly point out the main modifications.  The
primary change in the algorithm presented in this paper is the use of
target-specific expansions (Section~\ref{sec:ts}) to mediate the contributions
of the `directly evaluated' portion of the field---i.e., List~1, List~3~close,
and List~4~close. This replaces the formation of spherical harmonic expansions
at the QBX centers due to this portion of the field.  However, QBX local
expansions in spherical harmonics, formed at the QBX centers, remain part of the
algorithm as computational entities since they mediate the potential due to the
`far-field' as well as the field due to boxes in List~3~far.

There are two other less significant differences from the statement of the
algorithm in~\cite{gigaqbx3d}. First, we give a more precise statement of the
construction of the octree in this paper than that was given
in~\cite{gigaqbx3d}, when it comes to distinguishing what type of particles may
be `owned' by boxes. Second, we introduce the notion of a `List~3~far candidate
box,' and we modify the criteria for placement of boxes in List~3~far by introducing a source
count threshold.

\subsection{Notation}

In this section, we introduce the same notation as~\cite{gigaqbx3d} in support
of the precise statement of the modified algorithm. Let $b$ be a box in an
octree with center $c$. We use the notation $\closedbox(r, c)$ to denote the set
$\{ x \in \mathbb{R}^3 : \norminf{x - c} \leq r\}$ and $\closedball(r,c)$ to
denote the set $\{ x \in \mathbb{R}^3 : \norm{x - c} \leq r \}$.

We use $|b|$ to refer to the $\ell^\infty$ radius of the box, i.e.\ half the
box width.

The \emph{target confinement region (TCR)} of $b$, or $\tcr(b)$, is the set
$\closedball(\sqrt{3} |b| (1 + t_f), c)$, where $t_f \geq 0$ is the \emph{target
  confinement factor}.

The \emph{$k$-near~neighborhood} of $b$ is the region $\closedbox(|b|(1 + 2k),
c)$. The \emph{$k$-colleagues} of $b$ are same-level boxes contained inside the
$k$-near~neighborhood of $b$. In particular, $T_b$ denotes the set of
$2$-colleagues of a box $b$.

Two same-level boxes that are not $k$-colleagues are termed
\emph{$k$-well-separated}.

$\parent(b)$ denotes the parent of $b$.

The sets $\ancestors(b)$ and $\descendants(b)$ denote the sets of ancestors and
descendants of $b$. These are also defined for a set of boxes, as the union of
the boxes' sets of ancestors or descendants, respectively.

A box owning a point or QBX center target is called a \emph{target box}. A box
owning a source quadrature node is called a \emph{source box}. Ancestors of
target boxes are called \emph{target-ancestor boxes}.

Two boxes are \emph{adjacent} if the intersection of their boundaries is
non-empty, i.e.\ they share a common  face, edge, or corner.

We define a relation $\adequatesep$ over the set of boxes and target confinement
regions within the tree, with $a \adequatesep b$ to be read as `$a$ is
adequately separated from $b$, relative to the size of $a$'. We write $a
\adequatesep \tcr(b)$ for boxes $a$ and $b$ if the $\ell^2$ distance from the
center of $a$ to the boundary of $\tcr(b)$ is at least $3|a|$. We write $\tcr(a)
\adequatesep b$ for boxes $a$ and $b$ if the $\ell^2$ distance from the center
of $a$ to the boundary of $b$ is at least $3|a|(1 + t_f)$.
(As $\ell^2$ distance is bounded from below by $\ell^\infty$ distance, a
computationally convenient approximation for checking whether $\tcr(a) \adequatesep b$
is to check if the $\ell^\infty$ distance from the center of $a$ to the boundary
of $b$ is at least $3|a|(1+t_f)$.  This check is sufficient (but not necessary)
for $\tcr(a) \adequatesep b$. We use this approximation in
our implementation.)
We write $a \not \adequatesep b$ to denote the
negation of $a \adequatesep b$.

\subsection{Interaction Lists}%
\label{sec:ilists}

The computational domain of the algorithm is a box with equal-length sides (such
as the closed cube $[-1,1]^3$).
This `root box' is recursively partitioned into equal-sized octants, together forming an octree.
The root box contains all source quadrature nodes, QBX centers, and target points, which we
refer to generically as \emph{particles}, and also the entirety of each QBX
ball. Target points are classified as either\emph{conventional targets}
and \emph{QBX targets}, the latter requiring potential evaluation through a QBX
expansion. The class of \emph{thresholded particles} includes every particle type
except QBX targets. These are the particles to which the box particle count threshold
applies, as described below. Boxes, regardless of their having children, may `own' a subset of thresholded particles. QBX
targets are not owned by boxes but are instead `associated' with a QBX ball,
which in turn is `owned' by a box.

To construct the octree, boxes of the tree owning more than $\nmax > 0$ (a
user-set parameter) thresholded particles are iteratively subdivided,
transferring the ownership of particles into the child boxes, until the number
of thresholded particles per leaf (childless) box is below $\nmax$, or if all
potentially split boxes are empty due to constraints on QBX center placement. A
QBX center whose expansion ball cannot be contained in the TCR of the child box
is not transferred to the child and remains owned by the parent.

Information on the parts of the potential travels between boxes through
\emph{translation operators}, from sets of boxes indicated by \emph{interaction
lists}, which are lists of boxes attached to target or
target-ancestor boxes in the tree. These lists are based around a \emph{near-field}
$T_b$ consisting of the $2$-colleagues of a given box; i.e.\ the same-level nearest-neighbors
and second nearest-neighbors. For a more detailed description of the interaction
lists in~\algbrand, see~\cite{gigaqbx2d, gigaqbx3d}.

For a given box $b$, \emph{List 1} consists of interactions with adjacent
boxes.

\begin{definition}[List 1, $\ilist{1}{b}$~{\cite[Def.~2]{gigaqbx3d}}]\label{def:list-1}
  For a target box $b$, $\ilist{1}{b}$ consists of all leaf boxes from among
  $\descendants(b) \cup \{b\}$ and the set of boxes adjacent to $b$.
\end{definition}

\emph{List 2} consists of interactions with non-adjacent same-level boxes that
are descendants of the near-field of the parent.

\begin{definition}[List 2, $\ilist{2}{b}$~{\cite[Def.~3]{gigaqbx3d}}]\label{def:list-2}
  For a target or target-ancestor box $b$, $\ilist{2}{b}$ consists of the children of the
  $2$-colleagues of $b$'s parent that are $2$-well-separated from $b$.
\end{definition}

\emph{List 3} consists of interactions between non-adjacent boxes where the
source box is in the near-field of the target.

\begin{definition}[List 3, $\ilist{3}{b}$~{\cite[Def.~4]{gigaqbx3d}}]\label{def:list-3}
  For a target box $b$, a box $d \in \descendants(T_b)$ is in $\ilist{3}{b}$ if
  $d$ is not adjacent to $b$ and, for all $w \in \ancestors(d) \cap
  \descendants(T_b)$, $w$ is adjacent to $b$.
\end{definition}

\emph{List 4} consists of interactions between non-adjacent boxes where the
target box is in the near-field of the source.

\begin{definition}[List 4, $\ilist{4}{b}$,~{\cite[Def.~5]{gigaqbx3d}}]\label{def:list-4}
  For a target or target-ancestor box $b$, a source box $d$ is in List 4 of $b$ if
  $d$ is a $2$-colleague of some ancestor of $b$ and $d$ is adjacent to
  $\parent(b)$ but not $b$ itself. Additionally, a source box $d$ is in
  $\ilist{4}{b}$ if $d$ is a $2$-colleague of $b$ and $d$ is not adjacent to $b$.
\end{definition}

Lists 1--4 (cf.~\cite{carrier:1988:adaptive-fmm})
are generalizations of the interaction lists present in most variants of the FMM,
modified for a `2-away' near-field and the presence of targets in non-leaves.
To these definitions, we adjoin a set of `close' and `far' lists. The
purpose of these lists is to ensure that interactions directed at QBX centers
maintain sufficient separation so that intermediate translations involving QBX
centers have controlled accuracy. The field due to a `close' list is evaluated
directly (without the use of intermediate expansions), while the field due to a
`far' list is sufficiently separated to allow for the use of intermediate
expansions.

\emph{List 3 close} and \emph{List 3 far} consist of boxes from List 3 and their
descendants. A box is placed into one of these lists depending on whether it is
adequately separated from the TCR of the target box, and whether it exceeds a
certain `source count' threshold.

In order to define the close and far lists associated with $\ilist{3}{b}$,
we introduce the notion of a `List 3 far candidate' box.

\begin{definition}[List 3 far candidate box]
  For a target box $b$, a box $d$ is a \emph{List 3 far candidate of $b$} if
  \begin{inparaenum}[(a)]
    \item $d \in \ilist{3}{b} \cup \descendants(\ilist{3}{b})$,
    \item $d \adequatesep \tcr(b)$, and,
    \item for all $w \in \ancestors(d) \cap (\descendants(\ilist{3}{b}) \cup
      \ilist{3}{b})$, $w \not \adequatesep \tcr(b)$.
  \end{inparaenum}
\end{definition}

In other words, a List 3 far candidate box is a box in the near-field which is
the largest box which is adequately separated from the TCR of the target among
itself and its chain of ancestors. It follows that each ancestor chain of boxes
contains at most a single List 3 far candidate.  The `multipole threshold',
$\nmpole \geq 0$, a user-set threshold related to the cumulative number of
sources in the descendants of a List 3 far candidate, contributes to the placement of
candidates into List 3 far or List 3 close.

\begin{definition}[List 3 far, $\ilist{3far}{b}$, modified from~{\cite[Def.~7]{gigaqbx3d}}]\label{def:list-3-far}
  For a target box $b$, a List 3 far candidate $d$ of $b$ is said to be in
  $\ilist{3far}{b}$ if the cumulative number of sources owned by $d$ and its
  descendants is at least $\nmpole$.
\end{definition}

\begin{definition}[List 3 close, $\ilist{3close}{b}$, modified from~{\cite[Def.~6]{gigaqbx3d}}]\label{def:list-3-close}
  For a target box $b$, a leaf box $d$ is said to be in $\ilist{3close}{b}$ if
  $d \in \descendants(\ilist{3}{b}) \cup \ilist{3}{b}$ and one of the following
  is true:
  \begin{inparaenum}[(a)]
  \item $d \not \adequatesep \tcr(b)$,
  \item a box $w \in \{d\} \cup \ancestors(d)$ exists such that $w$ is a List 3 far
    candidate box of $b$ and the cumulative number of sources owned by $w$ and its
    descendants is less than $\nmpole$.
  \end{inparaenum}
\end{definition}

\emph{List 4 close} and \emph{List 4 far} consist of boxes that are in the List
4 of a target/target ancestor box or its chain of ancestors. A box is placed
into List 4 close if the TCR of the target/target ancestor box is not adequately
separated from it; otherwise it is placed in List 4 far.

\begin{definition}[List 4 close, $\ilist{4close}{b}$,~{\cite[Def.~8]{gigaqbx3d}}]\label{def:list-4-close}
  Let $b$ be a target or target-ancestor box. A box $d$ is in
  $\ilist{4close}{b}$ if for some $w \in \ancestors(b) \cup \{b\}$ we have $d
  \in \ilist{4}{w}$ and furthermore $\tcr(b) \not \adequatesep d$.
\end{definition}

\begin{definition}[List 4 far, $\ilist{4far}{b}$,~{\cite[Def.~9]{gigaqbx3d}}]\label{def:list-4-far}
  Let $b$ be a target or target-ancestor box. A box $d \in \ilist{4}{b}$ is in
  List 4 far if $\tcr(b) \adequatesep d$.  Furthermore, if $b$ has a parent, a
  box $d \in \ilist{4close}{\parent(b)}$ is in List 4 far if $\tcr(b)
  \adequatesep d$.
\end{definition}


\subsection{Formal Statement}%
\label{sec:formal-statement}

For completeness, we give a full statement of the~\algbrand\ algorithm using
target-specific expansions in Algorithm~\ref{alg:gigaqbx-ts}. We use the
following notation for `point' potentials (those formed without QBX mediation):
\begin{inparaenum}[(a)]
\item $\Potnear_b(t)$ denotes the potential at a target point $t$ due to all
  sources in $\ilist{1}{b} \cup \ilist{3close}{b} \cup \ilist{4close}{b}$, and
\item $\PotW{b}(t)$ denotes the potential at a target $t$ due to all sources in
  $\ilist{3far}{b}$.
\end{inparaenum}
For QBX-mediated potentials, we use the following notation, where $c$ is a QBX
center owned by a box $b$:
\begin{inparaenum}[(a)]
\item $\Lqbxnear{c}(t)$ denotes the (QBX) local expansion at a target $t$ due to
  all sources in $\ilist{1}{b} \cup \ilist{3close}{b} \cup \ilist{4close}{b}$,
\item $\LqbxW{c}(t)$ denotes the (QBX) local expansion at a target $t$ due to
  all sources in $\ilist{3far}{b}$, and
\item $\Lqbxfar{c}(t)$ denotes the (QBX) local expansion at a target $t$ due to
  sources not in $\ilist{1}{b} \cup \ilist{3}{b} \cup \ilist{4close}{b}$.
\end{inparaenum}

In Algorithm~\ref{alg:gigaqbx-ts}, the stages changed to use target-specific
expansions are indicated with a star ($\bigstar$).

\begin{algbreakable}{\algbrand~FMM with Target-Specific Expansions}%
  \label{alg:gigaqbx-ts}%
  \begin{algorithmic}
    \REQUIRE{The maximum number of thresholded
      FMM targets/sources $\nmax$ per box for
      octree refinement, a multipole threshold $\nmpole$, and a target
      confinement factor $t_f$ are chosen.}
    \REQUIRE{The input geometry and targets are preprocessed according
      to~\cite[Sec.~3]{gigaqbx3d}.}
    \REQUIRE{Based on the precision $\epsilon$ to be achieved, a QBX
      order $\pqbx$, an FMM order $\pfmm$, and an oversampled quadrature node
      count $\pquad$ are chosen.}
    \ENSURE{An accurate approximation to the potential at all target points
      is computed.}

    \COMMENT{$\bigstar$ indicates that a stage uses target-specific expansions.}

    \algstage{Stage 1: Build tree}
    {%
    \STATE{Create an octree on the computational domain containing all sources,
      targets, and QBX centers, as well as the entirety of each expansion ball.}
    \REPEAT{}
    \STATE{Subdivide each box owning more than $\nmax$ thresholded
      particles into eight children, pruning any empty child boxes. If a
      QBX center cannot be owned by the child box with target confinement
      factor $t_f$ due to its radius, it remains in the parent box.}
    \UNTIL{no box needs subdivision}
    }

    \algstage{Stage 2: Form multipoles}{%
    \FORALL{boxes $b$}
    \STATE{Form a $\pfmm$-th order multipole expansion $\Mpole_b$ centered at $b$ due to
      sources owned by $b$.}
    \ENDFOR{}
    \FORALL{boxes $b$ in postorder}
    \STATE{For each child of $b$, shift the center of the multipole expansion at
      the child to the center of $b$. Add the resulting expansions to $\Mpole_b$.}
    \ENDFOR{}
    }

    \algstage{$\bigstar$Stage 3: Evaluate direct interactions}{%
    \FORALL{boxes $b$}
    \STATE{For each non-QBX target $t$ owned by $b$, add to $\Potnear_b(t)$
      the contribution due to the interactions from sources owned by boxes in
      $\ilist{1}{b}$ to $t$.}
    \ENDFOR{}
    \FORALL{boxes $b$}
    \STATE{%
      For each QBX target $t$ associated to a QBX center $c$ owned by $b$,
      use target-specific expansions to add to $\Lqbxnear{c}(t)$ the contribution
      due to the interaction from all sources in $\ilist{1}{b}$.}
    \ENDFOR{}
    }

    \algstage{Stage 4: Translate multipoles to local expansions}{%
    \FORALL{boxes $b$}
    \STATE{For each box $d \in \ilist{2}{b}$, translate the multipole expansion
      $\Mpole_{d}$ to a local expansion centered at $b$. Add the resulting
      expansions to obtain $\Locfar_b$.}
    \ENDFOR{}
    }

    \algstage{$\bigstar$Stage 5(a): Evaluate direct interactions due to $\ilist{3close}{b}$}{%
    \STATE{Repeat Stage 3 with $\ilist{3close}{b}$ instead of $\ilist{1}{b}$.}
    }

    \algstage{Stage 5(b): Evaluate multipoles due to $\ilist{3far}{b}$}{%
    \FORALL{boxes $b$}
    \STATE{For each conventional target $t$ owned by $b$, evaluate the multipole
      expansion $\Mpole_{d}$ of each box $d \in \ilist{3far}{b}$
      to obtain $\PotW{b}(t)$.}
    \ENDFOR{}
    \FORALL{boxes $b$}
    \STATE{For each QBX center $c$ owned by $b$, compute the expansion
      $\LqbxW{c}$, due to
      the multipole expansion $\Mpole_{d}$ of each box $d \in \ilist{3far}{b}$.}
    \ENDFOR{}
    }

    \algstage{$\bigstar$Stage 6(a): Evaluate direct interactions due to $\ilist{4close}{b}$}{%
    \STATE{Repeat Stage 3 with $\ilist{4close}{b}$ instead of $\ilist{1}{b}$.}
    }

    \algstage{Stage 6(b): Form locals due to $\ilist{4far}{b}$}{%
    \FORALL{boxes $b$}
    \STATE{Convert the field of every source particle owned by boxes in $\ilist{4far}{b}$
      to a local expansion about $b$. Add to $\Locfar_b$.}
    \ENDFOR{}
    }

    \algstage{Stage 7: Propagate local expansions downward}{%
    \FORALL{boxes $b$ in preorder}
    \STATE{For each child $d$ of $b$, shift the center of the local expansions
      $\Locfar_b$ to the child. Add the resulting expansion to $\Locfar_d$.}
    \ENDFOR{}
    }

    \algstage{Stage 8: Form local expansions at QBX centers}{%
    \FORALL{boxes $b$}
    \STATE{For each QBX center $c$ owned by $b$, translate $\Locfar_b$ to
      $c$, obtaining $\Lqbxfar{c}$.}
    \ENDFOR{}
    }

    \algstage{Stage 9: Evaluate final potential at targets}{%
    \FORALL{boxes $b$}
    \STATE{For each non-QBX target $t$ owned by $b$, evaluate
      $\Locfar_b(t)$. Add $\Potnear_b(t)$, $\PotW{b}(t)$, and $\Locfar_b(t)$ to
      obtain the potential at $t$.}
    \ENDFOR{}
    \FORALL{boxes $b$}
    \STATE{For each QBX target $t$ associated to a QBX center $c$ owned by $b$,
      add $\Lqbxnear{c}(t)$, $\LqbxW{c}(t)$, and $\Lqbxfar{c}(t)$ to obtain the
      potential at $t$.}
    \ENDFOR{}
    }
  \end{algorithmic}
\end{algbreakable}

\FloatBarrier{}


\subsection{Accuracy and Scaling}

The accuracy bound for the error introduced in acceleration in
Algorithm~\ref{alg:gigaqbx-ts} is unchanged compared with the original
algorithm~\cite[Thm.~6]{gigaqbx3d}, due to the fact that target-specific
expansions are mathematically identical to their target-independent
counterparts. In essence, assuming the validity of a number of
hypotheses with strong numerical evidence~\cite[Hyp. 1--3]{gigaqbx3d},
the asymptotic acceleration error is
$O\left((3/4)^{\pfmm+1}\right)$ (the same as the `one-away' point FMM in three
dimensions~\cite{petersen_error_1995}) when $t_f \leq 0.84$.

One can analyze the scaling of Algorithm~\ref{alg:gigaqbx-ts} in a way that closely
parallels that of~\cite{gigaqbx3d}. The main difference in the analysis is connected
with the change in the asymptotic complexity due to the use of target-specific
expansions. In Table~\ref{tab:complexity-analysis}, we provide a summary of the
asymptotic complexity of the stages of the algorithm. The complexity is measured
in terms of asymptotic floating point operations (see Section~\ref{sec:perf}).

\def\nqbx{n_\mathrm{qbx}}
\def\nfmm{n_\mathrm{fmm}}

\def\EtoEcost{\nfmm^{3/2}}
\def\PtoEcost{\nfmm}
\def\PtoPTScost{\nqbx^{1/2}}
\def\PtoPcost{\nqbx}

\begin{DIFnomarkup}
\begin{table}
  \centering
  \caption{Complexity of each stage of Algorithm~\ref{alg:gigaqbx-ts}.}%
  \label{tab:complexity-analysis}
  \begin{tabular}{lll}
    \toprule
    Stage & Asymptotic Operation Count & Note \\
    \midrule
    Stage 1 & $O (NL)$ &  cf.~\cite[Tab.~1]{gigaqbx3d} \\
    Stage 2 & $O (N_S \PtoEcost + N_B \EtoEcost )$ & \ditto \\
    Stage 3 & $O ( (27 (N_C + N_S) \nmax + N_C M_C) \PtoPTScost )$ & Using TSQBX; cf.~\cite[Lem.~11]{gigaqbx3d} \\
    Stage 4 & $O (875 N_B \EtoEcost )$ & cf.~\cite[Lem.~12]{gigaqbx3d} \\
    Stage 5 & $O (N_C M_C \PtoPTScost + 124L N_S \nmax \PtoPTScost)$ & Using TSQBX; cf.~\cite[Lem.~13]{gigaqbx3d} \\
    Stage 6(a) & $O (250 N_C \nmax \PtoPTScost)$ & Using TSQBX; cf.~\cite[Lem.~15]{gigaqbx3d} \\
    Stage 6(b) &  $O (375 N_B \nmax \PtoEcost)$ & cf.~\cite[Lem.~15]{gigaqbx3d} \\
    Stage 7 & $O (8 N_B \EtoEcost)$ & cf.~\cite[Tab.~1]{gigaqbx3d} \\
    Stage 8 & $O (N_C \EtoEcost)$ & \ditto \\
    Stage 9 & $O (N_T \PtoPcost)$ & \ditto \\
    \bottomrule
  \end{tabular}
\end{table}
\end{DIFnomarkup}

We make a number of simplifying assumptions in the complexity model. The first
is that $\pqbx \leq \pfmm$, which is true in all practical situations as the error
introduced by FMM acceleration decreases more slowly with increasing order
compared with the error due to truncation of the QBX expansion
(cf.~Lemma~\ref{lem:qbx-truncation-3d}).
Secondly,
we assume the use of spherical harmonic expansions with `point-and-shoot'
translations (see Section~\ref{sec:perf}). Finally, we assume that the algorithm
is performing on-surface evaluation, in which the potential at every target is
mediated through a QBX center and there is at most one target per QBX center.
The presence of off-surface targets causes no significant changes to the complexity
analysis. We omit it for the sake of simplicity.

The following parameters are used as inputs to the complexity model. $N$, $N_T$,
$N_C$, $N_S$ refer respectively to the number of particles, targets, sources,
and centers. $N_B$ refers to the number of boxes and $L$ refers to the number of
levels in the tree. The quantity $M_C$, as found in~\cite[Lem.~8]{gigaqbx3d}, is
a measure of the average number of sources in a neighborhood of a QBX center,
where the size of the neighborhood is proportional to the ball size, defined
specifically as:
\[ M_C \coloneqq \frac{1}{N_C} \sum_{c \in C} \left| S \cap \closedbox\left(\frac{4 \sqrt{3} r_c}{t_f}, c\right) \right|, \]
where $C$ is the set of QBX centers, $S$ is the set of sources, and $r_c$ is the
expansion radius associated with the center $c$.  $\nfmm$ and $\nqbx$ refer to the number
of coefficients in an FMM and QBX expansion in spherical harmonics.

\subsection{Cost Considerations}%
\label{sec:params-affecting-performance}

While Table~\ref{tab:complexity-analysis} presents an understanding of the cost
dependence of the \algbrand~FMM's stages on the algorithmic parameters,
the cost analysis is asymptotic rather than a predictor of the time duration
from start to finish of the FMM, also known as the `wall time.' In this section,
we present a more qualitative understanding of the effect of algorithmic
parameters on cost. The observations in this section motivate an empirical study
of the cost of the \algbrand~FMM with target-specific expansions in the next
section.

When considering parameters that affect the cost of
Algorithm~\ref{alg:gigaqbx-ts}, a complication arises in that some algorithmic
parameters simultaneously have a large effect on both accuracy and
cost, so that modifying these parameters could result in output that is vastly different in accuracy.
A number of these parameters primarily affect the size of the QBX
near-field (which is covered by List 1, List 3 close, and List 3 far): $t_f$,
$\pquad$, and, to some extent, the input geometry itself. Additionally, the parameters
$\pqbx$ and $\pfmm$ do not affect the size of QBX near-field but nevertheless have
a major effect on accuracy. To ensure comparable levels of accuracy in the output
of the algorithm before and after optimization, we take the point of
view that parameters primarily affecting the accuracy of the layer potential
evaluation are fixed.

This leaves two parameters to be considered that primarily affect the cost of
the various stages of the algorithm, while leaving accuracy nearly
unchanged. The first of these parameters is $\nmax$, the maximum number of thresholded particles per
box. The main consideration in choosing a value of $\nmax$ is that as $\nmax$
increases, the number of boxes decreases, while simultaneously each box holds
more particles. The cost of those stages where the amount of total computational effort increases with the number of boxes,
such as Stage 4 (List 2), benefits from a reduction in the number of boxes.
On the other hand, stages
involving direct interactions require more work as $\nmax$ increases.

The second of these parameters is $\nmpole$. It was mentioned in~\cite{gigaqbx2d} as a degree
of freedom for optimization but not examined in detail. It is based on the observation that we
can avoid translating a multipole expansion into a QBX local expansion by
replacing it with direct interactions with the source particles whose field
makes up the multipole expansion. This is a less expensive evaluation strategy
for multipole expansions whose cumulative source count is below a `smallness' threshold.

\section{Experimental Results}%
\label{sec:results}

\subsection{Cost Model}%
\label{sec:perf}

In this section, we present a cost model for
Algorithm~\ref{alg:gigaqbx-ts} which estimates the amount of computational time used
by the algorithm in a
way that is reproducible while remaining predictive of actual machine computation time. To aid
the construction of a realistic and reproducible cost metric, the model
makes use of direct examination of the FMM tree and data
structures---a strategy we have found to yield data useful for a variety of
purposes. The model produces an
estimate of the total number of floating
point operations required for the algorithm on a particular input geometry.
Through the introduction of additional weight constants, we
use these counts to approximate the total amount of computational time used by
the algorithm, a quantity we refer to as the \emph{modeled process time}.
(The phrase `process time' in UNIX-type operating systems is used to
describe the total time spent executing process code, excluding the
time executing operating system code,
summed across all cores of a multi-core processor if relevant.)

The first step in obtaining an asymptotic estimate of the number of floating point
operations is to count the number of interactions of each category---e.g.,
local expansion formation,
multipole-to-local translation,
multipole evaluation,
etc.---performed by
the FMM by analyzing the interaction lists and tree. The model multiplies each
of these counts by a category-dependent symbolic expression, parametrized by the
number of coefficients in the FMM and QBX expansions, to obtain an asymptotic number
of floating point operations (e.g.,\ for forming a multipole/local expansion with $(p + 1)^2$ coefficients,
the asymptotic amount of work is modeled as $(p + 1)^2$).
Lastly, the
asymptotic number of floating point operations is multiplied by a
category-dependent `calibration constant', an empirically determined
parameter representing a ratio of running time in seconds to modeled floating
point operations.

The modeled process time per interaction category is
shown in Table~\ref{tab:model}.
To count the number of floating point operations,
the cost model
assumes the use of spherical harmonic expansions (see
Section~\ref{sec:sph-harm}). For each category of interaction, the modeled
number of floating point operations
is designed to be asymptotically correct to leading order.
We assume that translations between two expansions occur using a
`point-and-shoot' strategy --- this reduces the cost of the relevant
translations from $O((p+1)^4)$ to $O((p+1)^3)$ (see for
instance~\cite{gumerov:3d-translation-ops-comparison}) for homogeneous source
and target order $p$. The procedure for a `point-and-shoot' translation between
source order $p$ and target order $q$ is as follows:
\begin{enumerate}
\item At a cost of $O((p+1)^3)$, rotate the source expansion so that the translation
  direction is $z$-axis aligned.
\item At a cost of $O((p+1)(q+1)^2)$, translate the source expansion along the
  $z$-axis to the target expansion.
\item At a cost of $O((q+1)^3)$, rotate the target expansion back.
\end{enumerate}
Our model includes a term for each of these three stages. In the case of a
homogeneous source and target order, this entails using a leading factor of
$3$. Lastly, we model the cost of evaluation of a target-specific expansion as
$O(\pqbx + 1)$ floating-point operations, corresponding to the evaluation of
formula~\eqref{eqn:tsqbx} using recurrences for the Legendre polynomials.

The model is fitted to the results obtained from timing our implementation of
Algorithm~\ref{alg:gigaqbx-ts} on geometries of fixed QBX and FMM order. Timing
data is obtained by timing each stage of our implementation on a 20-core
2.30~GHz Intel Xeon~\mbox{E5-2650}~v3 machine. We use a least-squares fit to
obtain calibration coefficients from the timing data. Our implementation, which
uses double-precision floating point arithmetic throughout, is based on
FMMLIB~\cite{gimbutas_fmmlib} compiled by GCC~7.2.0 with \texttt{-Ofast} and
\texttt{-march=native} flags, and making use of shared memory parallelism via
OpenMP\@. The calibration coefficients obtained for the order pair $(\pqbx,
\pfmm) = (5, 15)$ are displayed in Table~\ref{tab:calibration-constants}. These
were obtained by fitting to the process times for the `urchin' geometries
$\gamma_3, \gamma_5$ (see Section~\ref{sec:urchin}).

A natural interpretation of these coefficients is that they represent the time
of a single `asymptotic flop' in their respective interaction category. It is therefore not
unreasonable to expect each flop to have an execution time roughly corresponding
to the inverse of the clock frequency of the processor. Since the processor we used
for our experiments has a clock frequency of 2.30~GHz, we expect and observe
asymptotic flop times of a magnitude around $10^{-9}$. Further, any major discrepancies
in their comparative magnitude may indicate a difference in implementation quality.
While most of the calibration
coefficients are of roughly the expected magnitude, the coefficient associated
with QBX-local-to-target evaluation is much larger due to inefficiencies our
implementation.  (Despite the high overhead of this evaluation, it does not
play a significant role in the overall cost of the scheme.)

Empirically, we have observed the model to give accurate cost estimates,
within a few percent of the true execution time, for the same QBX and FMM order
pair it is fitted to. As an example, Table~\ref{tab:perf-model-accuracy} gives
actual versus predicted process times for a sequence of `urchin' test geometries
introduction in Section~\ref{sec:urchin}, using the calibration constants from
Table~\ref{tab:calibration-constants}. This accuracy does not necessarily carry
over when differing values of $\pqbx$ and $\pfmm$ are used from those that the
model was fitted to, likely to the overhead of unmodeled lower-order costs in certain
interactions in our implementation. While this issue may
be addressed either extending the model to include more terms or by further
optimization work to reduce the lower-order costs in the implementation, in this
paper we handle this issue by re-fitting the model to each example that we use,
ensuring its fidelity as a predictor of modeled process time.

\def\cltop{c_\mathrm{l2p}}
\def\cmtop{c_\mathrm{m2p}}
\def\cmtom{c_\mathrm{m2m}}
\def\cptop{c_\mathrm{p2p}}
\def\cptom{c_\mathrm{p2m}}
\def\cptol{c_\mathrm{p2l}}
\def\cmtoqbxl{c_\mathrm{m2qbxl}}
\def\cptoqbxl{c_\mathrm{p2qbxl}}
\def\cltoqbxl{c_\mathrm{l2qbxl}}
\def\cqbxltop{c_\mathrm{qbxl2p}}
\def\cmtol{c_\mathrm{m2l}}
\def\cltol{c_\mathrm{l2l}}
\def\cptopts{c_\mathrm{ts}}

\begin{table}
  \centering
  \caption{Cost model used in this paper for evaluation of the
    scaling of Algorithm~\ref{alg:gigaqbx-ts}, where $\nqbx = (1+\pqbx)^2$ and $\nfmm=(1+\pfmm)^2$.}%
  \label{tab:model}
  \begin{tabular}{p{2in}l}
    \toprule
    Interaction & Modeled Process Time (s) \\
    \midrule
    Source $\to$ Local        & $\cptol \cdot \nfmm$ \\
    Source $\to$ Multipole    & $\cptom \cdot \nfmm$ \\
    Source $\to$ QBX Local    & $\cptoqbxl \cdot \nqbx$ \\
    Target-Specific QBX       & $\cptopts \cdot \nqbx^{1/2}$ \\
    Local  $\to$ Local        & $3 \cdot \cltol \cdot \nfmm^{3/2}$ \\
    Local  $\to$ QBX Local    & $\cltoqbxl \cdot \left(\nfmm^{3/2} + \nfmm^{1/2} \cdot \nqbx  + \nqbx^{3/2} \right)$ \\
    Multipole $\to$ Local     & $3 \cdot \cmtol \cdot \nfmm^{3/2}$ \\
    Multipole $\to$ Multipole & $3 \cdot \cmtom \cdot \nfmm^{3/2}$ \\
    Multipole $\to$ QBX Local & $\cmtoqbxl \cdot \left(\nfmm^{3/2} + \nfmm^{1/2} \cdot \nqbx  + \nqbx^{3/2} \right)$ \\
    QBX Local  $\to$ Target   & $\cqbxltop \cdot \nqbx$ \\
    \bottomrule
  \end{tabular}
  \vspace{1cm}
  \sisetup{table-format = 1.2e-1, exponent-product=\cdot}
  \caption{Calibration constants for the model in Table~\ref{tab:model},
    obtained for the order pair $(\pqbx, \pfmm) = (5, 15)$ on a particular
    machine using an implementation based on FMMLIB.}%
  \label{tab:calibration-constants}
  \begin{tabular}{lS}
    \toprule
    Constant & {Value} \\
    \midrule
    $\cptol$ & 1.10e-08 \\
    $\cptom$ & 1.24e-08 \\
    $\cptoqbxl$ & 1.42e-08 \\
    $\cptopts$ & 9.45e-09 \\
    $\cltol$ & 5.94e-09 \\
    \bottomrule
  \end{tabular}
  \quad
  \begin{tabular}{lS}
    \toprule
    Constant & {Value} \\
    \midrule
    $\cltoqbxl$ & 4.72e-09 \\
    $\cmtol$ & 3.24e-09 \\
    $\cmtom$ & 5.35e-09 \\
    $\cmtoqbxl$ & 3.37e-09 \\
    $\cqbxltop$ & 6.74e-07 \\
    \bottomrule
  \end{tabular}
  \vspace{1cm}
  \caption{Actual versus predicted process times using the model
    calibration constants in Table~\ref{tab:calibration-constants}
    for the `urchin' geometries in Section~\ref{sec:urchin},
    with fits obtained on the geometries $\gamma_3$, $\gamma_5$.}
  \label{tab:perf-model-accuracy}
  \begin{DIFnomarkup}
  \begin{tabular}{lrrrrr}
    \toprule
    & \multicolumn{5}{c}{Process Time (s)} \\
    \cmidrule(lr){2-6}
    \multicolumn{1}{c}{Kind} & \multicolumn{1}{c}{$\gamma_{2}$} & \multicolumn{1}{c}{$\gamma_{4}$} & \multicolumn{1}{c}{$\gamma_{6}$} & \multicolumn{1}{c}{$\gamma_{8}$} & \multicolumn{1}{c}{$\gamma_{10}$} \\
    \midrule
    Actual & 1457.07 & 5741.86 & 20705.79 & 46238.00 & 98943.60 \\
    Model & 1447.81 & 5731.83 & 20715.11 & 46261.86 & 98618.80 \\
    \bottomrule
  \end{tabular}
  \end{DIFnomarkup}
\end{table}

\subsection{Scaling and Balancing Study}

This section presents a numerical study of the impact of our optimizations
to the \algbrand\ algorithm. It is possible to use target-specific QBX as a direct
replacement for the version of the algorithm in~\cite{gigaqbx3d}, leaving all
other parameters in the algorithm unchanged, but the
simple adaptation of a few algorithmic parameters has the potential to
drastically affect performance characteristics by 
changing the proportion of work performed by different translation operators
in the FMM\@.
In the experiments of this section,
we apply this procedure to
Algorithm~\ref{alg:gigaqbx-ts} in order improve algorithmic cost,
by reducing performance bottlenecks, in a process known as `balancing'.

For all experiments, we report on (a) the total modeled process
time as well as (b) the modeled process times of the algorithmic stages that involve the
interaction lists (see Section~\ref{sec:ilists}). In the remainder of this
section, we freely use the name of an interaction list to refer to the stage it
is associated with.

\subsubsection{Test Geometry}%
\label{sec:urchin}

\begin{figure}
  \centering
  \includegraphics{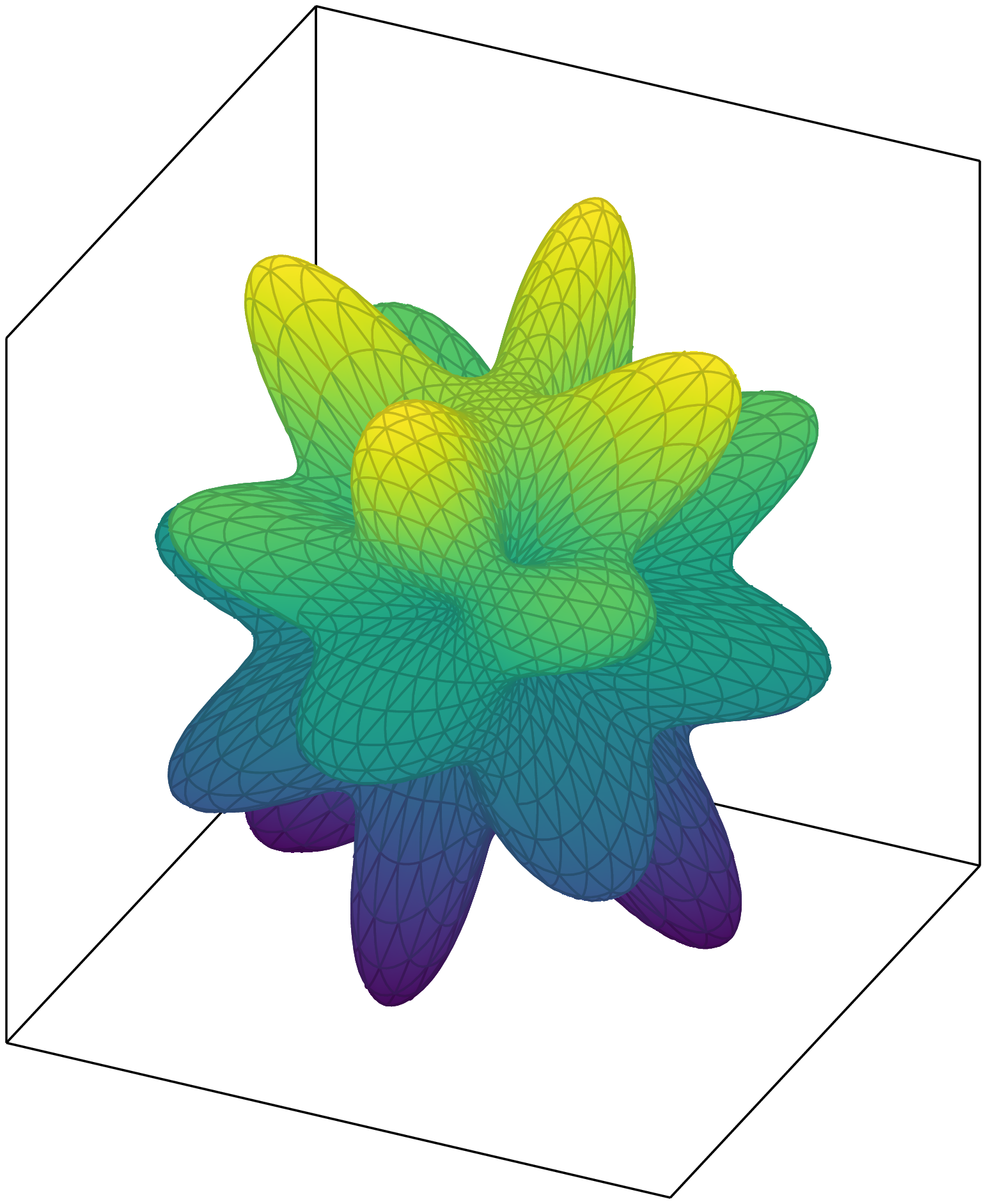}
  \caption{A mesh of the `urchin' test geometry $\gamma_8$.}%
  \label{fig:urchin}
  \vspace{1cm}
  \includegraphics{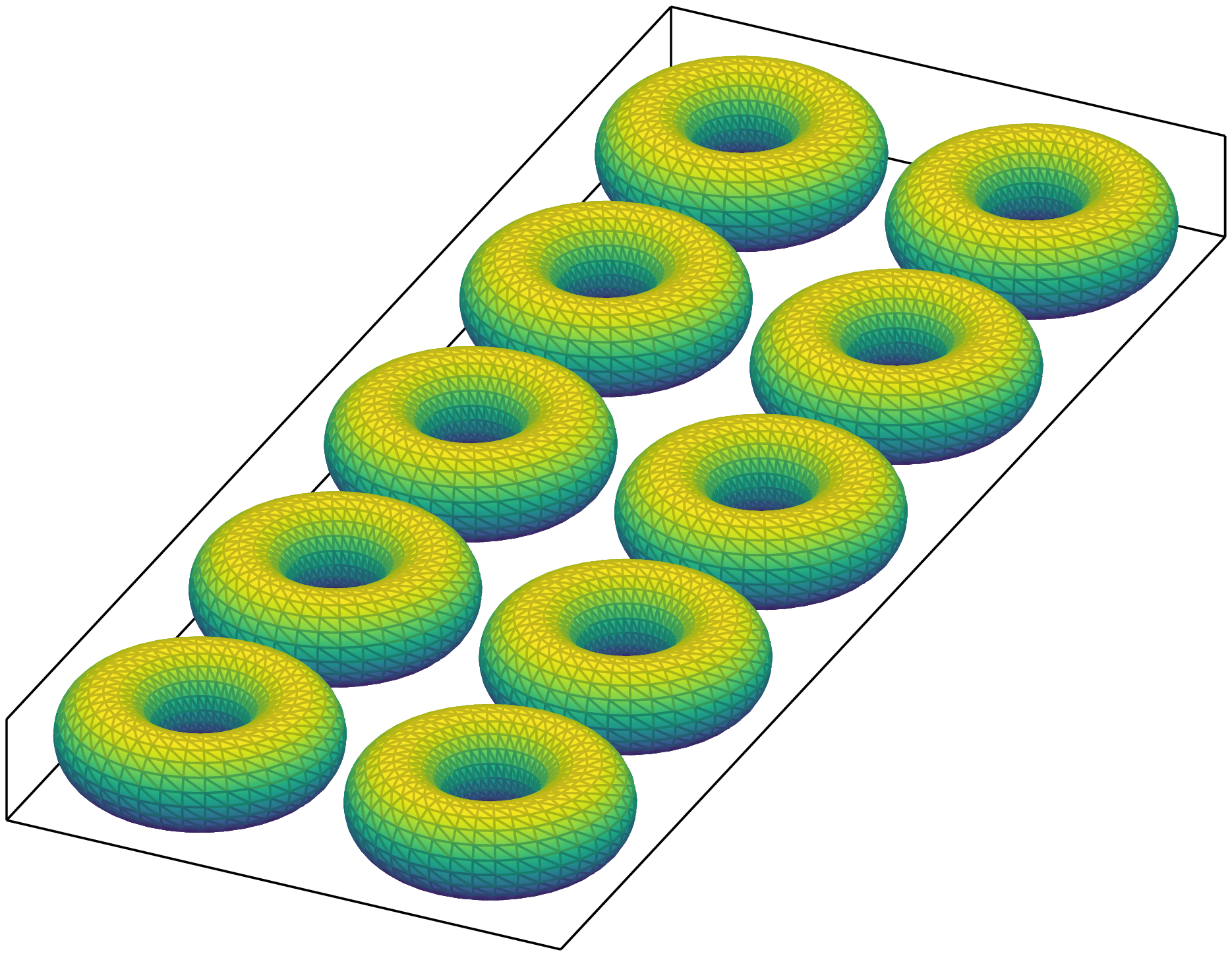}
  \caption{A mesh of the `torus grid' test geometry $\tau_{10}$.}%
  \label{fig:donuts}
\end{figure}

\def\urchinorder{8}
\def\minparticles{\num{4116770}}
\def\maxparticles{\num{178529030}}
\def\urchinquadnodes{295}

As a first test of the scaling behavior and cost of our algorithm, we use a family of test
geometries $\gamma_k$ that we call `urchins', parametrized by $k \in
\mathbb{N}$. These surfaces are given in spherical coordinates as $(r_k, \theta,
\phi)$, using the definition of spherical harmonics in~\eqref{eqn:sph-harm},
letting $r_k$ vary with $\theta$ and $\phi$ as follows:
\begin{align*}
  \label{eq:urchin-warping}
  r_k(\theta,\phi) &= 0.2 + \frac{\Re Y_{k}^{\lfloor k/2\rfloor}(\theta,\phi) - m_k}{M_k-m_k},\\
  M_k &= \max_{\theta\in[0,\pi],\phi\in[0,2\pi]} \Re Y_{k}^{\lfloor k/2\rfloor}(\theta,\phi),\\
  m_k &= \min_{\theta\in[0,\pi],\phi\in[0,2\pi]} \Re Y_{k}^{\lfloor k/2\rfloor}(\theta,\phi).
\end{align*}
Figure~\ref{fig:urchin} gives a visual impression of the geometry
$\gamma_8$. The mesh of the geometry consists of triangular elements whose
parameter mapping function is an \urchinorder{}th degree polynomial. The
construction starts with the image of an icosahedral mesh under the mapping
function $r_k$. To ensure that the piecewise polynomial elements accurately
represent the geometry, an iterative refinement procedure is applied to the
mapped elements, and at each iteration any refined elements are nodally
reevaluated. The details of this procedure can be found
in~\cite[Sec.~6]{gigaqbx3d}.

In our numerical experiments, we use the family of geometries $\gamma_2,
\gamma_4, \ldots, \gamma_{10}$, which range in size from
\minparticles--\maxparticles\ particles, where a `particle' is a source, target,
or QBX center. We use $(\pqbx, \pfmm) = (5, 15)$. As reported
in~\cite{gigaqbx3d}, this corresponds to about five digits of accuracy for
evaluating Green's identity~\cite[Thm.~6.5]{kress:2014:integral-equations},
which we use as a proxy test for accuracy in layer potential evaluation.  The upsampled
quadrature rule has \urchinquadnodes~nodes per element, and we use $t_f = 0.9$.

\subsubsection{Cost Evaluation}%
\label{sec:perf-evaluation}

\begin{figure}
  \centering
  \input{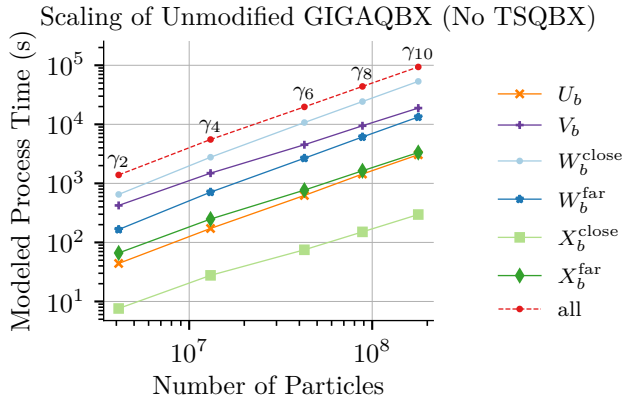}
  \vspace*{-3mm}
  \caption{Scaling of the unmodified \algbrand~algorithm (with spherical harmonic expansions instead of target-specific expansions)
    for evaluation of the Laplace single-layer potential on the sequence of
    `urchin' geometries $\gamma_2, \gamma_4, \ldots, \gamma_{10}$.}%
  \label{fig:complexity-urchin}
\end{figure}

\def\optnmaxnots{96}
\def\optnmpolenots{40}
\def\optnmaxts{992}
\def\optnmpolets{280}

In this section, we consider the cost of evaluation of the on-surface value of
the Laplace single-layer potential on the `urchin' geometry family. We establish
a `baseline' cost by modeling version of the \algbrand~algorithm as described
in~\cite{gigaqbx3d} using the same framework from Section~\ref{sec:perf}.
The results are shown in
Figure~\ref{fig:complexity-urchin}~(cf.~\cite[Fig.~12]{gigaqbx3d}, which reports
modeled floating point operations for the same geometries). The total reported
time (under `all') includes the time contribution from all stages of the
algorithm. These results are obtained using the model in
Section~\ref{sec:perf}. To obtain this set of results, we chose a set of
balancing parameters that minimized modeled process time as described in
Section~\ref{sec:perf}, which are $\nmax = \optnmaxnots$ and $\nmpole =
\optnmpolenots$.

\begin{figure}
  \centering
  \input{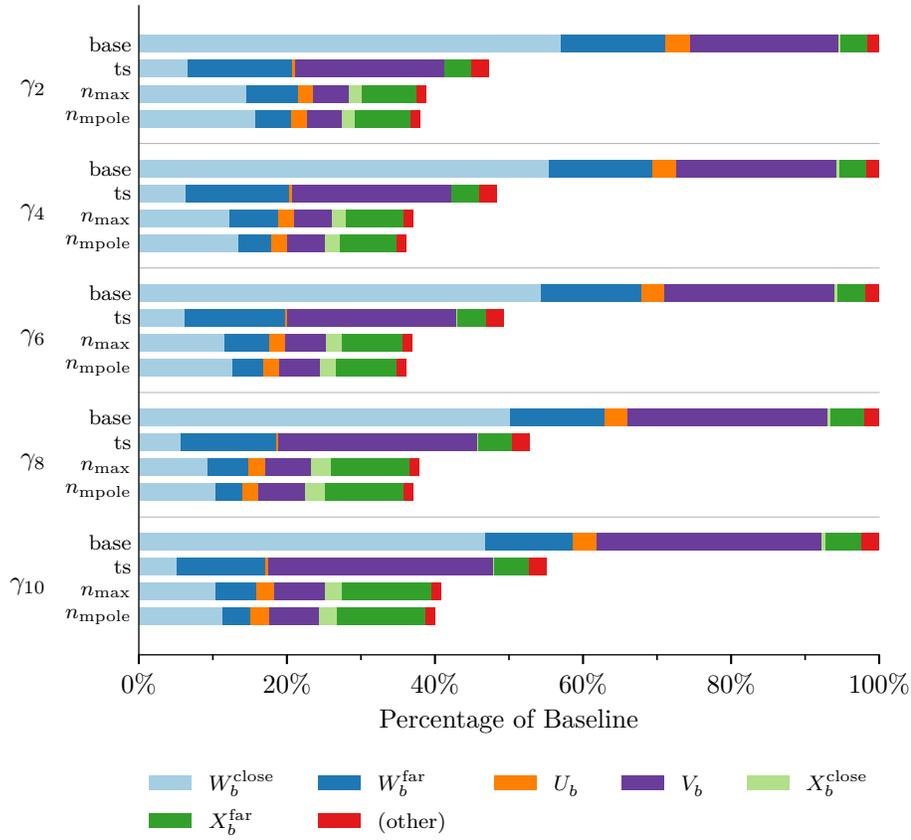}
  \caption{Cumulative impact of a sequence of optimizations, applied to the
    evaluation of the Laplace single-layer potential for the `urchin' family of
    geometries, on modeled process time as well as time for individual
    stages. `base' denotes the baseline time; `ts' denotes the result of using
    target-specific expansions; `$\nmax$' denotes the result of rebalancing
    $\nmax$; and `$\nmpole$' denotes the result of rebalancing $\nmpole$.}%
  \label{fig:urchin-opts}
\end{figure}

To assess the impact of each of the optimizations in this paper, we present
results of three stages of cumulative optimizations. These are summarized in
Figure~\ref{fig:urchin-opts}.

\def\tsrawspeedup{9}

\paragraph*{Using TSQBX Without Rebalancing.} The first optimization we consider is making use of
target-specific QBX expansions, leaving all other algorithmic parameters
constant. According to Section~\ref{sec:perf}, assuming an average of one target
per center, this should give a speedup of
\[ \cptoqbxl \left(1+\pqbx\right)^2 / \left(\cptopts (1+\pqbx)\right) \]
in the evaluation of List 3 close, List 4 close, and List 1, while leaving the
cost of the other stages of the algorithm unchanged. This turns out to be the
case. The data in Table~\ref{tab:calibration-constants} predict a cost reduction
by a factor of about~\tsrawspeedup\ according to this formula. In
Figure~\ref{fig:urchin-opts}, a cost reduction of this magnitude is evident in
$\ilist{1}{b}$, $\ilist{3close}{b}$, and $\ilist{4close}{b}$ under the label
`ts'. On average, the overall time is 50\% of the baseline.

\begin{figure}[ht]
  \centering
  \input{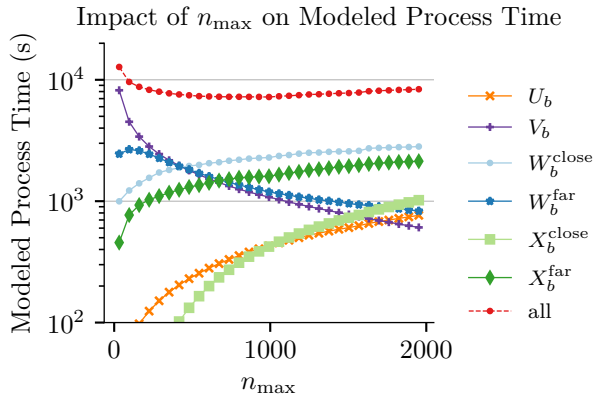}%
  \vspace*{-5mm}
  \caption{Modeled process time of the \algbrand~algorithm and its various stages,
    versus $\nmax$, for the geometry $\gamma_6$, using target-specific QBX\@.}%
  \label{fig:nmax-urchin-study}
\end{figure}

\paragraph*{Rebalancing $\bm{\nmax}$.} Next, we consider the effect of varying
$\nmax$. Using $\gamma_6$ as a reference geometry, we conduct a study measuring
total modeled process time against $\nmax$. The results are presented in
Figure~\ref{fig:nmax-urchin-study}. Based on these results we choose $\nmax =
\optnmaxts$, which appears to be empirically near-optimal. The effect on
the running time of the various stages is given under `$\nmax$' in
Figure~\ref{fig:urchin-opts}.  On average, the running time is 38\% of the
baseline.

As we can see, the increase in $\nmax$ has a complex effect on the time for the
different stages of the algorithm. Perhaps most easily explained, the time
associated with List 2 decreases. This is consistent with a decrease in the
overall number of boxes in the tree that may be expected with an increase in
$\nmax$.

We also observe an increase in the amount of work for List 4 far. The amount of
work for List 4 far is proportional to both $\nmax$ and the number of boxes
(cf.~\cite[Lem.~15]{gigaqbx3d},~Table~\ref{tab:complexity-analysis}).
While the number of boxes decreases as $\nmax$
increases, the increase in $\nmax$ appears to have the dominant effect on the cost
of this stage.

The proportion of time associated with target-specific expansions appears to increase
roughly proportionally with $\nmax$. This is plausible as the cost estimates from
previous work (cf.~\cite[Lem.~11,13,15]{gigaqbx3d},~Table~\ref{tab:complexity-analysis})
show that at least a
portion of the time in these stages is proportional to $\nmax$.

It is more difficult to explain why the time associated List 3 far decreases.
However, as the decrease in the time associated with List 3 far closely
mirrors that of List 2 (see Figure~\ref{fig:nmax-urchin-study}), it is likely to
be related to the decrease in the number of boxes.

\paragraph*{Rebalancing $\bm{\nmpole}$.} The last optimization we consider is the
balancing of the constant $\nmpole$.  Using an experiment similar to that
of the one done to rebalance $\nmax$, we obtain a threshold value of $\nmpole =
\optnmpolets$.

The cost model can also be used to obtain this threshold. From the
cost model, we expect approximately minimal cost when
\[
    \nmpole \approx \frac{\cmtoqbxl}{\cptopts}
    \frac{
      (1+\pfmm)^3 + (1+\pqbx)^2(1+\pfmm) + (1+\pqbx)^3
    }{
      (1+ \pqbx)
    }.
\]
Based on the data from Table~\ref{tab:calibration-constants}, this value
is about 291, so that the theoretical value and the
empirically obtained value from the cost model are nearly in agreement.

The impact of increasing $\nmpole$ to reduce the amount of work done in
$\ilist{3far}{b}$, shifting the corresponding work to $\ilist{3close}{b}$, for a
reduction in the total modeled process time time to 37\% of the original, on average.
This can be seen in
Figure~\ref{fig:urchin-opts} under the label `$\nmpole$'.

\paragraph*{Remarks.} Overall, the cumulative impact of these optimizations is to
reduce the overall running time to an average of 37\% of the original. As the
percentage reduction in cost is nearly uniform across geometries,
the scaling characteristics in geometry size are essentially unchanged
from the baseline version in Figure~\ref{fig:complexity-urchin}, in that
we observe an approximately linear scaling of the algorithm cost with the
particle count.

From the relative costs in Figure~\ref{fig:urchin-opts}, it is evident that this
version of the algorithm is more evenly `balanced' than the
baseline, in that the proportion of time spent in individual stages is distributed
more uniformly. Specifically, the cost contribution of List 3 close is now
significantly closer to the other interaction lists. The dominant costs appear
to be List 3 close and List 4 far.

\subsubsection{An Example with Higher Accuracy and Simpler Geometry}%
\label{sec:donut}

\begin{figure}
  \centering
  \input{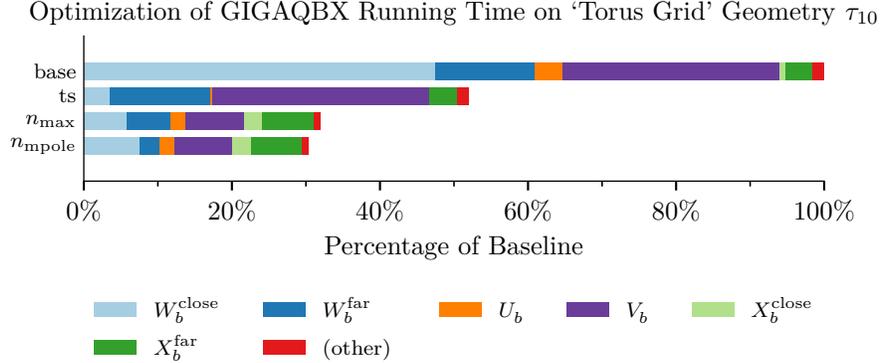}
  \caption{Cumulative impact of a sequence of optimizations, applied to the
    evaluation of the Laplace single-layer potential for the `torus grid'
    geometry $\tau_{10}$, on modeled process time as well as time for individual
    stages.}%
  \label{fig:donut-opts}
\end{figure}

To demonstrate the generality of the optimizations mentioned in this paper, we
repeat the cost evaluation in Section~\ref{sec:perf-evaluation} on a different, simpler geometry
while choosing parameters to attain higher accuracy. The geometry $\tau_1$ is a torus parametrized as the image of the set
$(u, v) \in [0, 2\pi) \times [0, 2\pi)$ under the mapping
\begin{align*}
  \label{eqn:torus}
  x &= \cos u \left(1 + 2 \cos v \right) \\
  y &= \cos u \left(1 + 2 \cos v \right) \\
  z &= 2 \sin v.
\end{align*}
The `torus grid' geometries $\tau_{2k}$, $k \in \mathbb{N}$, are obtained by
spacing $2k$ copies of $\tau_1$ on a $2 \times k$ grid, with a uniform spacing
of 0.6. See Figure~\ref{fig:donuts} for a visual impression of $\tau_{10}$. To
obtain the initial mesh of the torus, we tile the parameter domain into $40
\times 20$ contiguous rectangles, and then subdivide each rectangle into two
triangles. We represent each triangle as the image of an triangular reference element under a polynomial mapping of degree $8$,
with 295
quadrature nodes per element, and use $t_f = 0.9$.

We focus on the evaluation of the Laplace single-layer potential on the geometry
$\tau_{10}$. We use $(\pqbx, \pfmm) = (9, 20)$.  From the initial mesh for
$\tau_{10}$, refinement produces a mesh with \num{16000} stage-1 elements and
\num{64000} stage-2 elements, for a total number of about 19 million source
particles. Since the effects of closeness to touching would reflect in the count
of stage-1 elements, the number of stage-1 elements ($40 \cdot 20 \cdot 2 \cdot
10$) indicates that refinement due to closeness to touching of the different
components is not necessary. A test of Green's identity with the QBX and FMM
order parameters yields about eight digits of accuracy.

We optimize for evaluation on $\tau_{10}$. The modeled process time for the
baseline \algbrand~FMM, without using target-specific expansions, is minimized
with a choice of $\nmax = 96$ and $\nmpole = 40$. The baseline modeled process
time is \num{18968} seconds. We apply the sequence of optimizations mentioned in
Section~\ref{sec:perf-evaluation}. The results are shown in
Figure~\ref{fig:donut-opts}. The use of target-specific expansions reduces the
modeled process time to 52\% of the baseline, which is shown under the label
`ts'. Using an empirically determined value of $\nmax = 928$ reduces the time to
32\% of the baseline, which is shown under the label `$\nmax$'. Finally, using an
empirically determined value of $\nmpole = 420$ reduces the time to 30\% of the
baseline, which is shown under the label `$\nmpole$'.

In many ways, the cost characteristics of this example are similar to
that of the examples in the previous section, although, due to the effects of
higher order FMM and QBX expansions, the cost reduction is more
significant. Like for the urchin, the dominant cost for the baseline example is in
the near-field evaluations. After rebalancing, the dominant costs are List 3
close and List 4 far.

\subsection{A Large-Scale BVP for the Helmholtz Equation}

\begin{figure}
  \centering
  \includegraphics[width=\textwidth]{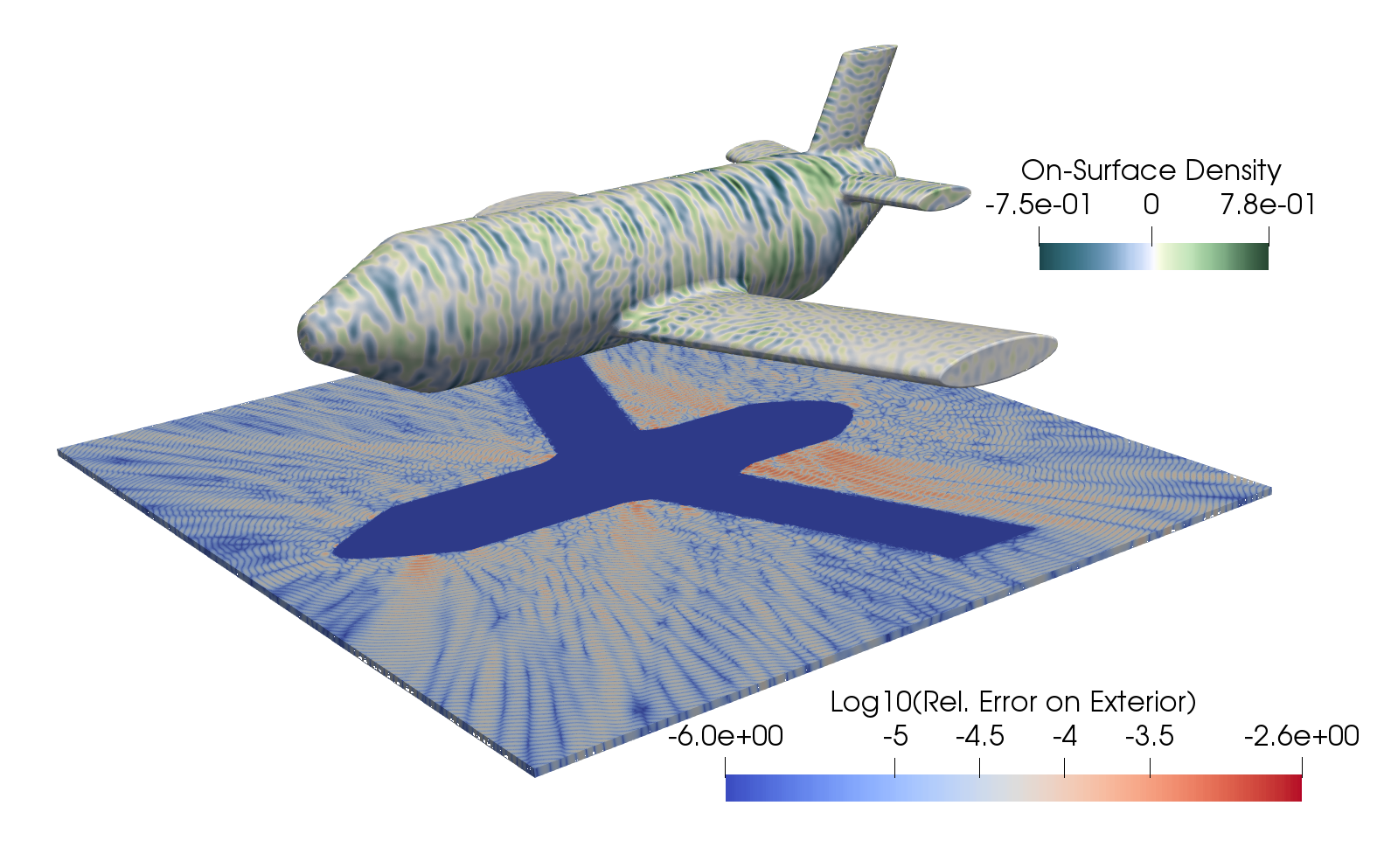}
  \caption{Visualization of the solution to an exterior Dirichlet problem for
    the Helmholtz equation solved on the `plane' geometry. The surface geometry
    is shaded according to the real part of the solved-for density $\mu$. The slice
    positioned below the geometry in the visualization is taken from the level
    of the `wing' of the plane and is shaded according to the logarithm of the
    observed relative error in the exterior of the volume. The maximum observed
    relative error anywhere in the exterior is about $2.5 \cdot 10^{-3}$.}%
  \label{fig:betterplane}
  \vspace{1cm}
  \input{timing-plane.pgf}
  \caption{Cumulative impact of a sequence of optimizations, applied to the
    evaluation of the Helmholtz double-layer potential, $k = 20$, for the `plane' geometry,
    on modeled process time as well as time for individual stages.}%
  \label{fig:plane-opts}
\end{figure}

As a final demonstration of the broad applicability of the optimizations in this
paper, we present a numerical example involving a large-scale boundary value
problem with complex geometry. Thus far, we have only discussed the use of
target-specific expansions for the Laplace equation, but everything we have
stated in this paper has an analogue involving the Helmholtz kernel (see
Appendix~\ref{sec:ts-expansions}). To demonstrate this, we solve the exterior
Dirichlet problem for the Helmholtz equation
\begin{DIFnomarkup}
\begin{alignat*}{2}
  \left( \triangle + k^2 \right) u(x) &= 0 & \quad & x \in \mathbb R^3 \setminus \Omega, \\
  u(x)                                &= f(x) & \quad & x \in \Gamma, \\
  \lim_{\norm{x} \to\infty} \norm{x} \left(\frac{\partial}{\partial \norm{x}} - ik \right) u(x)
                                      &= 0,
\end{alignat*}
\end{DIFnomarkup}
where $\Omega \subset \mathbb R^3$ is a closed, bounded region with smooth
boundary $\Gamma = \partial \Omega$. The solution $u$ to the boundary value
problem uses a Brakhage-Werner representation~\cite{brakhage_uber_1965}
\[ u \coloneqq i \alpha \mathcal{S} \mu - \mathcal{D} \mu, \]
where $\alpha \in \mathbb{R}$ and the double-layer operator $\mathcal{D}$ for a
Green's function $\kernel$ is defined as
\[
  (\mathcal{D} \mu)(x) \coloneqq \int_{\Gamma}
  \frac{\partial \kernel(x,y)}{\partial
   \nu(y)} \mu(y) \, dS(y).
\]

Our choice of geometry $\Omega$ is derived from
\texttt{surface-3d/betterplane.brep} from~\cite{geometries_repo_2018}; this is
the same source geometry as the example in~\cite[Sec.~6.2]{gigaqbx3d}. The
surface mesh of the geometry consists of triangular elements of degree two
mapping functions obtained with Gmsh~\cite{geuzaine_gmsh_2009}. The stage-1
discretization consists of \num{60638} elements and the stage-2 discretization
consists of \num{91526} elements. We use 150 quadrature nodes per element, for a total
of about 14 million source points, and $t_f = 0.9$.

The geometry, visualized in Figure~\ref{fig:betterplane}, has a bounding box of
size approximately $19 \times 20 \times 6.5$ units. We choose a Helmholtz
parameter of $k = 20$. The reference solution is obtained using a point
potential of sources places near the `tail' of the geometry. We modify the
quadrature scheme to use $L^2$-weighted degrees of freedom as a way to improve
conditioning of the discrete problem~\cite{bremer_nystrom_2011}.

Using the \algbrand~FMM with target-specific expansions and near-optimal parameter
values to drive the matrix-vector products in GMRES~\cite{saad_gmres_1986}, our method
attains a decrease in the residual norm by a factor of $10^{-6}$ in 120
iterations. We use $\pqbx = 4$ and we choose the FMM order to ensure a relative
FMM error below $10^{-5}$. This entails a minimum FMM order of $\pfmm = 13$,
but, as typical in an implementation of the Helmholtz
FMM~\cite{gumerov-helmholtz-fmm}, also requires increasing $\pfmm$ with the box
size to maintain accuracy tolerances.  We use Helmholtz translation operators
from FMMLIB3D~\cite{gimbutas_fmmlib,gimbutas_fast_2009}.  According to
comparison with the reference solution using point evaluations in the exterior
of the volume, the scheme obtains a relative $\ell^\infty$ error of about $2.5
\cdot 10^{-3}$.

We apply the optimization sequence described in Section~\ref{sec:perf-evaluation}
to obtain algorithmic parameters to minimize the modeled process time for the
double-layer evaluation with target-specific expansions. The latter operator is
the dominant cost of a GMRES iteration.

The parameter values that obtain near-optimal modeled cost for the `baseline'
version without target-specific expansions are $\nmax = 100$ and $\nmpole =
40$. The baseline modeled process time of the double-layer operator is \num{28900} seconds.
Figure~\ref{fig:plane-opts} shows the effect of the various
optimizations on the modeled process time for the double-layer operator.
Using target-specific expansions reduces the
modeled process time to 73\% of the baseline, shown under the label
`ts'. An empirically determined value of $\nmax = 300$ reduces the time to
60\% of the baseline, shown under the label `$\nmax$'. Lastly, an
empirically determined value of $\nmpole = 150$ reduces the time to 59\% of the
baseline, shown under the label `$\nmpole$'. The cost decreases in this example
are smaller compared with the examples of~Sections~\ref{sec:perf-evaluation}~and~\ref{sec:donut},
since the QBX expansion order is lower and the cost of direct evaluations is not as
significant in the baseline.

We also measured the `wall time', or the actual duration from start to finish of the
operator, using the baseline and optimized versions. The baseline wall time of the
double-layer operator is about \num{2300} seconds, and the wall time for the single layer
is about \num{1500} seconds.
The sequence of optimizations above reduces the two times to about \num{1600} seconds and
\num{900} seconds respectively. Hence, the use of target-specific expansions reduces
the solve time for this problem to 66\% of the baseline.
Compared with process time, the wall time is not reduced
as much because not every stage of our implementation is efficiently parallelized. This
remains as future work.

\section{Conclusion}%
\label{sec:conclusion}

This paper examines optimizations to reduce the cost of the \algbrand~FMM in
three dimensions. The main cost, as reported in~\cite{gigaqbx3d} and verified in
this paper, is typically due to the increased number of near-field interactions
as compared with the `point' version of the FMM\@. To reduce the impact of the
cost of the near-field interactions, we consider an acceleration strategy using target-specific expansions.
We develop a version of the \algbrand~FMM that uses
target-specific expansions for near-field evaluation. We also demonstrate a
cost model for the \algbrand~FMM that accurately estimates the total process time
used on a shared-memory system. This model aids in an efficient assessment of
the impact of various choices of algorithmic parameters affecting
computational cost. We find that, in order to make the best use of target-specific expansions,
algorithmic parameters should be modified to effectively shift more of the work
onto target-specific near-field interactions. We demonstrate that this strategy
produces speedups of 1.7--3.3$\times$ in practice on examples involving the
Laplace and Helmholtz kernels using a variety of QBX and FMM orders.






A key strength of target-specific expansions is their general
applicability. However, the requisite modification to the algorithm, including
the mathematical derivation of target-specific expansions for specific kernels
and the optimization of algorithmic parameters as done in this paper, requires
manual intervention. A technical question remains as to what extent these
tasks can be automated. One current subject of our investigation is the use of
techniques from symbolic computing and/or numerical linear algebra to develop
automated ways to generate optimal complexity target-specific expansions for
various kernels. A second subject is a better combination of automation with the
cost model in this paper in order to enable efficient selection of
parameters, such as the maximum number of particles per box, for arbitrary input geometries.

\section*{Acknowledgments}

The authors' research was supported by the National Science Foundation under
grants DMS-1418961 and DMS-1654756 as well as the Department of Computer Science
at the University of Illinois at Urbana-Champaign.  Any opinions, findings, and
conclusions, or recommendations expressed in this article are those of the
authors and do not necessarily reflect the views of the National Science
Foundation; NSF has not approved or endorsed its content.  Part of the work was
performed while the authors were participating in the HKUST-ICERM workshop
`Integral Equation Methods, Fast Algorithms and Their Applications to Fluid
Dynamics and Materials Science' held in 2017. The authors would also
like to thank the anonymous reviewer for helpful comments.

\appendix
\section{Target-Specific Expansions}%
\label{sec:ts-expansions}

For completeness, we describe how to obtain target-specific expansions for the
kernels of the single-layer potential $\mathcal{S}\mu$, its normal derivative
$\mathcal{S}'\mu$, and the double-layer potential $\mathcal{D}\mu$ for a Green's
function $\mathcal{K}$, where $\mathcal{K}$ is the Laplace or Helmholtz kernel:
\begin{align*}
  (\mathcal{S} \mu)(x) &= \int_{\Gamma} \kernel(x, y) \mu(y) \, dS(y), \\
  (\mathcal{S}' \mu)(x) &= \int_{\Gamma} \frac{\partial \kernel(x,y)}{\partial
    \nu(x)} \mu(y) \, dS(y), \\
  (\mathcal{D} \mu)(x) &= \int_{\Gamma} \frac{\partial \kernel(x,y)}{\partial
    \nu(y)} \mu(y) \, dS(y).
\end{align*}
Consider a source $s \in \mathbb{R}^3$, center $c \in \mathbb{R}^3$, and target
$t \in \mathbb{R}^3$ satisfying $\norm{t - c} \leq \norm{s - c}$. Let $\gamma$
be the angle between $s - c$ and $t - c$ (Figure~\ref{fig:tsqbx}).

In the case of the Laplace kernel
$\mathcal{G}$~\eqref{eqn:laplace-green-function}, the $p$-th order
target-specific expansion $\mathcal{G}^{(p)}$ takes the form
\begin{equation}
  \mathcal{G}^{(p)}(t,s) = \frac{1}{4 \pi}
  \sum_{n=0}^p \frac{\norm{t-c}^n}{\norm{s-c}^{n+1}} P_n(\cos \gamma).
\end{equation}
To obtain the target-specific expansion of the normal derivative of the single-layer
kernel, we require the gradient with respect to the target $t$. For the
double-layer kernel, the gradient should be taken with respect to the source
$s$. Recalling that $\cos \gamma = (t - c)^\top (s - c) / (\norm{t-c}
\norm{s-c})$, we use the chain rule to obtain
\begin{align}
  \begin{aligned}
    \nabla_t \mathcal{G}^{(p)}(t,s) &=
    \frac{1}{4 \pi} \sum_{n=1}^p \frac{\norm{t-c}^{n-1}}{\norm{s-c}^{n+1}}
    \left(
    n\frac{t-c}{\norm{t-c}} P_n(\cos \gamma)
    \right. \\
    &\qquad \quad \left.
    +\left[
      \frac{s-c}{\norm{s-c}} -
      \frac{t-c}{\norm{t-c}} \cos \gamma
      \right] P_n'(\cos \gamma)
    \right),
  \end{aligned}\\
  \begin{aligned}
    \nabla_s \mathcal{G}^{(p)}(t,s) &=
    \frac{1}{4 \pi} \sum_{n=0}^p \frac{\norm{t-c}^n}{\norm{s-c}^{n+2}}
    \left(
    -(n+1)\frac{s-c}{\norm{s-c}} P_n(\cos \gamma)
    \right. \\
    &\qquad \quad \left.
    +\left[
      \frac{t-c}{\norm{t-c}} -
      \frac{s-c}{\norm{s-c}} \cos \gamma
      \right] P_n'(\cos \gamma)
    \right).
  \end{aligned}
\end{align}

In the Helmholtz case, the kernel is
\[
  \mathcal{G}_k(t, s) = \frac{e^{ik\norm{t - s}}}{4 \pi \norm{t - s}}.
\]
For $|t| < |s|$, the addition for the Helmholtz kernel takes the form (see
for instance~\cite[eq.~(10.60.1,10.60.2)]{nist:dlmf}):
\[
  \mathcal{G}_k(t, s) =
  \frac{ik}{4 \pi} \sum_{n=0}^\infty (2n + 1) j_n(k\norm{t}) h_n(k\norm{s})
  P_n(\cos \theta),
\]
where the functions $j_n$ and $h_n$ are spherical Bessel and Hankel functions of
the first kind (see~\cite[eq.~(10.47.3)]{nist:dlmf}
and~\cite[eq.~(10.47.5)]{nist:dlmf}). This implies that the $p$-th order target-specific
expansion of the Helmholtz kernel takes the form
\begin{equation}
  \mathcal{G}^{(p)}_k(t, s) = \frac{ik}{4 \pi}
  \sum_{n=0}^p (2n + 1) j_n(k\norm{t-c}) h_n(k\norm{s-c}) P_n(\cos \gamma).
\end{equation}
The gradients with respect to the target $t$ and the source $s$ are
\begin{align}
  \begin{aligned}
    \nabla_t \mathcal{G}_k^{(p)}(t,s) &=
    \frac{ik}{4 \pi}
    \sum_{n=0}^p
    (2n+1)\frac{h_n(k\norm{s-c})}{\norm{t-c}}
    \left( \vphantom{\left[\frac{t-c}{\norm{t-c}}\right]}
    k (t - c) j_n'(k\norm{t-c}) P_n(\cos \gamma)
    \right. \\
    &\qquad \qquad \left.
    + \left[
      \frac{s-c}{\norm{s-c}} -
      \frac{t-c}{\norm{t-c}} \cos \gamma
      \right] j_n(k\norm{t-c}) P_n'(\cos \gamma)
    \right),
  \end{aligned}\\
  \begin{aligned}
    \nabla_s \mathcal{G}_k^{(p)}(t,s) &=
    \frac{ik}{4 \pi}
    \sum_{n=0}^p
    (2n+1)\frac{j_n(k\norm{t-c})}{\norm{s-c}}
    \left( \vphantom{\left[\frac{t-c}{\norm{t-c}}\right]}
    k (s - c) h_n'(k\norm{s-c}) P_n(\cos \gamma)
    \right. \\
    &\qquad \qquad \left.
    + \left[
      \frac{t-c}{\norm{t-c}} -
      \frac{s-c}{\norm{s-c}} \cos \gamma
      \right] h_n(k\norm{s-c}) P_n'(\cos \gamma)
    \right).
  \end{aligned}
\end{align}

Recurrences are available for computing the functions $P_n$, $j_n$, $h_n$, and
their derivatives rapidly~\cite[Sec.~10.51,~14.10]{nist:dlmf}.

\section{Software and Reproducibility}%
\label{sec:software}

The results in this paper were obtained using
\textbf{Pytential}~\cite{pytential}, a Python package for the
evaluation of layer potentials in two and three dimensions and the
solution of related integral equations. We used revision
\verb|54360f5| (tagged \verb|tsqbx|), available at
\url{https://github.com/inducer/pytential}.

In addition, we have prepared a Docker image from which all experiments in
this paper can be automatically reproduced with included software and scripts,
available at \url{https://doi.org/10.5281/zenodo.3523410}. The
code for the experiments is available at \url{https://doi.org/10.5281/zenodo.3542253}.

\printbibliography{}

@article{carrier:1988:adaptive-fmm,
  author       = {Carrier, J. and Greengard, L. and Rokhlin, V.},
  doi          = {10.1137/0909044},
  fjournal     = {Society for Industrial and Applied Mathematics. Journal on Scientific and Statistical Computing},
  % issn         = {0196-5204},
  journaltitle = {SIAM J. Sci. Statist. Comput.},
  mrclass      = {65C20 (70-08 78-08 78A35 82-08)},
  mrnumber     = {945931},
  number       = {4},
  pages        = {669--686},
  title        = {A fast adaptive multipole algorithm for particle simulations},
  % url          = {http://dx.doi.org/10.1137/0909044},
  volume       = {9},
  year         = {1988},
}

@article{epstein:2013:qbx-error-est,
  author       = {Epstein, Charles L. and Greengard, Leslie and Klöckner, Andreas},
  doi          = {10.1137/120902859},
  fjournal     = {SIAM Journal on Numerical Analysis},
  % issn         = {0036-1429},
  journaltitle = {SIAM J. Numer. Anal.},
  mrclass      = {65R20 (31B10 65N38 65N80)},
  mrnumber     = {3106484},
  mrreviewer   = {Michael J. Carley},
  number       = {5},
  pages        = {2660--2679},
  title        = {On the convergence of local expansions of layer potentials},
  % url          = {http://dx.doi.org/10.1137/120902859},
  volume       = {51},
  year         = {2013},
}

@article{klockner:2013:qbx,
  author       = {Klöckner, Andreas and Barnett, Alexander H. and Greengard, Leslie and O'Neil, Michael},
  doi          = {10.1016/j.jcp.2013.06.027},
  fjournal     = {Journal of Computational Physics},
  % issn         = {0021-9991},
  journaltitle = {J. Comput. Phys.},
  mrclass      = {65D30 (45Exx)},
  mrnumber     = {3101510},
  pages        = {332--349},
  title        = {Quadrature by expansion: a new method for the evaluation of layer potentials},
  % url          = {http://dx.doi.org/10.1016/j.jcp.2013.06.027},
  volume       = {252},
  year         = {2013},
}

@article{rachh:2017:qbx-fmm,
  author       = {Rachh, Manas and Klöckner, Andreas and O'Neil, Michael},
  doi          = {10.1016/j.jcp.2017.04.062},
  fjournal     = {Journal of Computational Physics},
  % issn         = {0021-9991},
  journaltitle = {J. Comput. Phys.},
  mrclass      = {Prelim},
  mrnumber     = {3667635},
  pages        = {706--731},
  title        = {Fast algorithms for {Q}uadrature by {E}xpansion {I}: {G}lobally valid expansions},
  % url          = {http://dx.doi.org/10.1016/j.jcp.2017.04.062},
  volume       = {345},
  year         = {2017},
}

@article{rahimian:2017:qbkix,
  author="Rahimian, Abtin and Barnett, Alex and Zorin, Denis",
  title="Ubiquitous evaluation of layer potentials using Quadrature by Kernel-Independent Expansion",
  journal="BIT Numerical Mathematics",
  year="2017",
  month="11",
  day="06",
  abstract="We introduce a quadrature scheme---QBKIX ---for the ubiquitous high-order accurate evaluation of singular layer potentials associated with general elliptic PDEs, i.e., a scheme that yields high accuracy at all distances to the domain boundary as well as on the boundary itself. Relying solely on point evaluations of the underlying kernel, our scheme is essentially PDE-independent; in particular, no analytic expansion nor addition theorem is required. Moreover, it applies to boundary integrals with singular, weakly singular, and hypersingular kernels. Our work builds upon quadrature by expansion, which approximates the potential by an analytic expansion in the neighborhood of each expansion center. In contrast, we use a sum of fundamental solutions lying on a ring enclosing the neighborhood, and solve a small dense linear system for their coefficients to match the potential on a smaller concentric ring. We test the new method with Laplace, Helmholtz, Yukawa, Stokes, and Navier (elastostatic) kernels in two dimensions (2D) using adaptive, panel-based boundary quadratures on smooth and corner domains. Advantages of the algorithm include its relative simplicity of implementation, immediate extension to new kernels, dimension-independence (allowing simple generalization to 3D), and compatibility with fast algorithms such as the kernel-independent FMM.",
  % issn="1572-9125",
  doi="10.1007/s10543-017-0689-2",
  %url="https://doi.org/10.1007/s10543-017-0689-2"
}

@article{afklinteberg:2016:quadrature-est,
  author       = {af Klinteberg, Ludvig and Tornberg, Anna-Karin},
  doi          = {10.1007/s10444-016-9484-x},
  fjournal     = {Advances in Computational Mathematics},
  % issn         = {1019-7168},
  journaltitle = {Adv. Comput. Math.},
  mrclass      = {65D30},
  mrnumber     = {3598839},
  mrreviewer   = {Alexandru Ioan Mitrea},
  number       = {1},
  pages        = {195--234},
  title        = {Error estimation for quadrature by expansion in layer potential evaluation},
  % url          = {http://dx.doi.org/10.1007/s10444-016-9484-x},
  volume       = {43},
  year         = {2017},
}

@book{kress:2014:integral-equations,
  author    = {Kress, Rainer},
  doi       = {10.1007/978-1-4614-9593-2},
  edition   = {Third edition},
  isbn      = {978-1-4614-9593-2},
  mrclass   = {45A05 (45-02 45L05 47G10 65R20)},
  mrnumber  = {3184286},
  pages     = {xvi+412},
  publisher = {Springer, New York},
  series    = {Applied Mathematical Sciences},
  title     = {Linear integral equations},
  % url       = {http://dx.doi.org/10.1007/978-1-4614-9593-2},
  volume    = {82},
  year      = {2014},
}

@article{pouransari:2015:adaptive-fractal-fmm,
  author       = {Pouransari, Hadi and Darve, Eric},
  fjournal     = {SIAM Journal on Scientific Computing},
  % issn         = {1064-8275},
  doi          = {10.1137/140962681},
  journaltitle = {SIAM J. Sci. Comput.},
  mrclass      = {65R10 (28A80 70F10)},
  mrnumber     = {3340201},
  mrreviewer   = {Peter Junghanns},
  number       = {2},
  pages        = {A1040--A1066},
  title        = {Optimizing the adaptive fast multipole method for fractal sets},
  % url          = { https://doi.org/10.1137/140962681},
  volume       = {37},
  year         = {2015},
}

@article{siegel2017local,
  title = "A local target specific quadrature by expansion method for evaluation of layer potentials in 3D",
  journal = "Journal of Computational Physics",
  volume = "364",
  pages = "365 - 392",
  year = "2018",
  % issn = "0021-9991",
  doi = "10.1016/j.jcp.2018.03.006",
  % url = "http://www.sciencedirect.com/science/article/pii/S002199911830158X",
  author = "Michael Siegel and Anna-Karin Tornberg",
  keywords = "Layer potentials, Integral equations, Quadrature by expansion, Exterior Dirichlet problem, Spherical harmonics expansions, Multiply-connected domain"
}

@article{saad_gmres_1986,
  title = "{GMRES: A Generalized Minimal Residual Algorithm for Solving Nonsymmetric Linear Systems}",
  volume = {7},
  doi = {10.1137/0907058},
  number = {3},
  journal = {{SIAM} Journal on Scientific and Statistical Computing},
  author = {Saad, Youcef and Schultz, Martin H.},
  month = jul,
  year = {1986},
  pages = {856--869},
}

@article{gigaqbx2d,
  title = "A fast algorithm with error bounds for Quadrature by Expansion",
  journal = "Journal of Computational Physics",
  volume = "374",
  pages = "135 - 162",
  year = "2018",
  doi = "10.1016/j.jcp.2018.05.006",
  author = "Matt Wala and Andreas Klöckner"
}

@article{gigaqbx3d,
  title = "A fast algorithm for Quadrature by Expansion in three dimensions",
  journal = "Journal of Computational Physics",
  volume = "388",
  pages = "655 - 689",
  year = "2019",
  % issn = "0021-9991",
  doi = "10.1016/j.jcp.2019.03.024",
  %url = "http://www.sciencedirect.com/science/article/pii/S0021999119302074",
  author = "Matt Wala and Andreas Klöckner",
  keywords = "Fast algorithms, Fast multipole method, Integral equations, Quadrature, Singular integrals, Three dimensional problems",
}

@techreport{gumerov:3d-translation-ops-comparison,
  author={Gumerov, Nail A. and Duraiswami, Ramani},
  title="{Comparison of the efficiency of translation operators used in the fast multipole method for the 3D Laplace equation}",
  institution="{University of Maryland Institute of Advanced Computer Studies}",
  year={2005},
}

@article{xiao_numerical_2010,
  title = {A numerical algorithm for the construction of efficient quadrature rules in two and higher dimensions},
  volume = {59},
  % issn = {08981221},
  doi = {10.1016/j.camwa.2009.10.027},
  number = {2},
  journal = {Computers \& Mathematics with Applications},
  author = {Xiao, Hong and Gimbutas, Zydrunas},
  month = jan,
  year = {2010},
  pages = {663--676},
  file = {ScienceDirect - Computers & Mathematics with Applications \: A numerical algorithm for the construction of efficient quadrature rules in two and higher dimensions:/home/andreas/Zotero/storage/8HNTDAT9/science.html:text/html}
}

@article{vioreanu_spectra_2014,
  title = {Spectra of {Multiplication} {Operators} as a {Numerical} {Tool}},
  volume = {36},
  % issn = {1064-8275},
  doi = {10.1137/110860082},
  abstract = {We introduce a numerical procedure for the construction of interpolation and quadrature formulae on bounded convex regions in the plane. The construction is based on the behavior of spectra of certain multiplication operators and leads to nodes which are inside a prescribed convex region in \$\{{\textbackslash}mathbb R\}{\textasciicircum}2\$. The resulting interpolation schemes are numerically stable and the quadrature formulae have positive weights and almost (but not quite) optimal numbers of nodes. The performance of the algorithm is illustrated by several numerical examples.},
  number = {1},
  % urldate = {2018-05-09},
  journal = {SIAM Journal on Scientific Computing},
  author = {Vioreanu, B. and Rokhlin, V.},
  month = jan,
  year = {2014},
  pages = {A267--A288},
  file = {Full Text PDF:/home/andreas/Zotero/storage/ZK4JMMD2/Vioreanu und Rokhlin - 2014 - Spectra of Multiplication Operators as a Numerical.pdf:application/pdf;Snapshot:/home/andreas/Zotero/storage/T7JP3GX3/110860082.html:text/html}
}

@article{brakhage_uber_1965,
  title = {Über das {Dirichletsche} {Außenraumproblem} für die {Helmholtzsche} {Schwingungsgleichung}},
  volume = {16},
  % issn = {0003-889X},
  doi = {10.1007/BF01220037},
  number = {1},
  journal = {Archiv der Mathematik},
  author = {Brakhage, Helmut and Werner, Peter},
  year = {1965},
  keywords = {Mathematics and Statistics},
  pages = {325--329},
  file = {SpringerLink Full Text PDF:/home/andreas/data/zotero/storage/VT7TJTTT/Brakhage und Werner - 1965 - Über das Dirichletsche Außenraumproblem für die He.pdf:application/pdf;SpringerLink Snapshot:/home/andreas/data/zotero/storage/UKIPT3AT/abstract.html:text/html}
}

@article{bremer_nystrom_2011,
  title = {On the {Nyström} discretization of integral equations on planar curves with corners},
  journal = {Applied and Computational Harmonic Analysis},
  author = {Bremer, J.},
  year = {2011},
  doi={10.1016/j.acha.2011.03.002}
}

@article{geuzaine_gmsh_2009,
  author = {Christophe Geuzaine and Jean‐François Remacle},
  title = {Gmsh: A 3‐D finite element mesh generator with built‐in pre‐ and post‐processing facilities},
  journal = {International Journal for Numerical Methods in Engineering},
  year = {2009},
  volume = {79},
  number = {11},
  pages = {1309-1331},
  doi = {10.1002/nme.2579}
}

@misc{geometries_repo_2018,
  author={Andreas Klöckner},
  title={A Repository of Sample Geometries},
  url={https://github.com/inducer/geometries/},
  note={Retrieved at revision \texttt{a869fc3ad}}
}

@misc{nist:dlmf,
  key = "{\relax DLMF}",
  title = "{\it NIST Digital Library of Mathematical Functions}",
  howpublished = "https://dlmf.nist.gov/, Release 1.0.18 of 2018-03-27",
  url = "https://dlmf.nist.gov/",
  note = "F.~W.~J. Olver, A.~B. {Olde Daalhuis}, D.~W. Lozier, B.~I. Schneider,
         R.~F. Boisvert, C.~W. Clark, B.~R. Miller and B.~V. Saunders, eds."
}

@article{petersen_error_1995,
  title = {Error estimates for the fast multipole method. {II}. {The} three-dimensional case},
  volume = {448},
  copyright = {Scanned images copyright © 2017, Royal Society},
  % issn = {0962-8444, 2053-9177},
  % url = {http://rspa.royalsocietypublishing.org/content/448/1934/401},
  doi = {10.1098/rspa.1995.0024},
  abstract = {The fast multipole method is a new and interesting method for computing long range interactions in particle systems. Although the method already has been implemented and compared to traditional and other methods with respect to accuracy and speed, no accurate error estimates of the fast multipole method have been given. In this paper we develop an explicit though complicated form for the error in the three dimensional case, and derive an estimate for it. The estimate has a simple analytic form which will allow its use in tuning the method for best efficiency and for comparison of the method with other methods at the same accuracy.},
  language = {en},
  number = {1934},
  % urldate = {2018-05-13},
  journal = {Proc. R. Soc. Lond. A},
  author = {Petersen, Henrik G. and Smith, E. R. and Soelvason, D.},
  month = mar,
  year = {1995},
  pages = {401--418},
  file = {Snapshot:/home/andreas/Zotero/storage/6VXV4QGM/401.html:text/html}
}

@article{gimbutas_fast_2009,
  title = {A fast and stable method for rotating spherical harmonic expansions},
  volume = {228},
  % issn = {0021-9991},
  % url = {http://www.sciencedirect.com/science/article/pii/S0021999109002691},
  doi = {10.1016/j.jcp.2009.05.014},
  abstract = {In this paper, we present a simple and efficient method for rotating a spherical harmonic expansion. This is a well-studied problem, arising in classical scattering theory, quantum mechanics and numerical analysis, usually addressed through the explicit construction of the Wigner rotation matrices. We show that rotation can be carried out easily and stably through “pseudospectral” projection, without ever constructing the matrix entries themselves. Existing fast algorithms, based on recurrence relations, are subject to a variety of instabilities, limiting the effectiveness of the approach for expansions of high degree.},
  number = {16},
  %urldate = {2012-07-07},
  journal = {Journal of Computational Physics},
  author = {Gimbutas, Z. and Greengard, L.},
  month = sep,
  year = {2009},
  keywords = {Rotation matrix, Spherical harmonics, Fast multipole method},
  pages = {5621--5627},
  file = {ScienceDirect Full Text PDF:/home/andreas/data/zotero/storage/XKBE9SUK/Gimbutas und Greengard - 2009 - A fast and stable method for rotating spherical ha.pdf:application/pdf;ScienceDirect Snapshot:/home/andreas/data/zotero/storage/2HQFZ5XV/S0021999109002691.html:text/html}
}

@misc{gimbutas_fmmlib,
  title={FMMLIB3D},
  author={Gimbutas, Zydrunas and Greengard, Leslie},
  url={https://github.com/zgimbutas/fmmlib3d},
  note={Retrieved at revision \texttt{339e93bbc}}
}

@book{kellogg,
  author    = "Oliver D. Kellogg",
  title     = "Foundations of Potential Theory",
  publisher = "Dover Publications",
  year      = "2000",
  address   = "New York, NY"
}

@book{gumerov-helmholtz-fmm,
  title = {Fast Multipole Methods for the Helmholtz Equation in Three Dimensions},
  author = {Nail A. Gumerov and Ramani Duraiswami},
  address = {Amsterdam},
  publisher = {Elsevier Science},
  year = {2004},
  doi = {10.1016/B978-0-08-044371-3.X5000-5},
  series = {Elsevier Series in Electromagnetism},
}

@inproceedings{Choi:2014:CGH:2588768.2576787,
 author = {Choi, Jee and Chandramowlishwaran, Aparna and Madduri, Kamesh and Vuduc, Richard},
 title = {A {CPU-GPU} Hybrid Implementation and Model-Driven Scheduling of the Fast Multipole Method},
 booktitle = {Proceedings of Workshop on General Purpose Processing Using GPUs},
 series = {GPGPU-7},
 year = {2014},
 isbn = {978-1-4503-2766-4},
 location = {Salt Lake City, UT, USA},
 pages = {64:64--64:71},
 articleno = {64},
 numpages = {8},
 doi = {10.1145/2576779.2576787},
 acmid = {2576787},
 publisher = {ACM},
 address = {New York, NY, USA},
}

@inproceedings{Chandramowlishwaran:2012:BAT:2312005.2312039,
 author = {Chandramowlishwaran, Aparna and Choi, JeeWhan and Madduri, Kamesh and Vuduc, Richard},
 title = {Brief Announcement: Towards a Communication Optimal Fast Multipole Method and Its Implications at Exascale},
 booktitle = {Proceedings of the Twenty-fourth Annual ACM Symposium on Parallelism in Algorithms and Architectures},
 series = {SPAA '12},
 year = {2012},
 isbn = {978-1-4503-1213-4},
 location = {Pittsburgh, Pennsylvania, USA},
 pages = {182--184},
 numpages = {3},
 doi = {10.1145/2312005.2312039},
 acmid = {2312039},
 publisher = {ACM},
 address = {New York, NY, USA},
}

@techreport{agullo:hal-01474556,
  TITLE = {{Modeling Irregular Kernels of Task-based codes: Illustration with the Fast Multipole Method}},
  AUTHOR = {Agullo, Emmanuel and Bramas, B{\'e}renger and Coulaud, Olivier and Stanisic, Luka and Thibault, Samuel},
  %URL = {https://hal.inria.fr/hal-01474556},
  TYPE = {Research Report},
  NUMBER = {RR-9036},
  PAGES = {35},
  INSTITUTION = {{INRIA Bordeaux}},
  YEAR = {2017},
  MONTH = Feb,
  KEYWORDS = {Mathematical Software ;  Modeling and simulation ;  Parallel computing methodologies ;   fast multipole method ;  runtime system ;  task-based programming},
  PDF = {https://hal.inria.fr/hal-01474556/file/rapport.pdf},
  HAL_ID = {hal-01474556},
  HAL_VERSION = {v1},
}

@misc{pytential,
    title={Pytential: a software package for the evaluation of layer potentials},
    author={Klöckner, Andreas and Wala, Matt},
    url={https://github.com/inducer/pytential},
    note={Accessed October 2018}
}

@article {cheng:1999:3d-fmm,
    AUTHOR = {Cheng, H. and Greengard, L. and Rokhlin, V.},
     TITLE = {A fast adaptive multipole algorithm in three dimensions},
   JOURNAL = {J. Comput. Phys.},
  FJOURNAL = {Journal of Computational Physics},
    VOLUME = {155},
      YEAR = {1999},
    NUMBER = {2},
     PAGES = {468--498},
      % ISSN = {0021-9991},
   MRCLASS = {65N55},
  MRNUMBER = {1723309},
MRREVIEWER = {So-Hsiang Chou},
       DOI = {10.1006/jcph.1999.6355},
       %URL = {https://doi.org/10.1006/jcph.1999.6355},
}

@article{anderson:1992:fmm,
    AUTHOR = {Anderson, Christopher R.},
     TITLE = {An implementation of the fast multipole method without
              multipoles},
   JOURNAL = {SIAM J. Sci. Statist. Comput.},
  FJOURNAL = {Society for Industrial and Applied Mathematics. Journal on
              Scientific and Statistical Computing},
    VOLUME = {13},
      YEAR = {1992},
    NUMBER = {4},
     PAGES = {923--947},
      % ISSN = {0196-5204},
   MRCLASS = {65C99 (34B27 35J05)},
  MRNUMBER = {1166169},
       DOI = {10.1137/0913055},
       %URL = {https://doi.org/10.1137/0913055},
}

@article {singer:1995:fmm,
    AUTHOR = {Singer, J. K.},
     TITLE = {Parallel implementation of the fast multipole method with
              periodic boundary conditions},
   JOURNAL = {East-West J. Numer. Math.},
  FJOURNAL = {East-West Journal of Numerical Mathematics},
    VOLUME = {3},
      YEAR = {1995},
    NUMBER = {3},
     PAGES = {199--216},
      % ISSN = {0928-0200},
   MRCLASS = {65Y05 (65N99 70F10)},
  MRNUMBER = {1359410},
MRREVIEWER = {N. Gass},
}

@article {zhao:1991:fmm,
    AUTHOR = {Zhao, Feng and Johnsson, S. Lennart},
     TITLE = {The parallel multipole method on the {C}onnection {M}achine},
   JOURNAL = {SIAM J. Sci. Statist. Comput.},
  FJOURNAL = {Society for Industrial and Applied Mathematics. Journal on
              Scientific and Statistical Computing},
    VOLUME = {12},
      YEAR = {1991},
    NUMBER = {6},
     PAGES = {1420--1437},
      % ISSN = {0196-5204},
   MRCLASS = {65Y05 (65C20 70-08 70F10)},
  MRNUMBER = {1129654},
       DOI = {10.1137/0912077},
       %URL = {https://doi.org/10.1137/0912077},
}

@techreport{hu:1996:accuracy,
  title={On the Accuracy of Poisson's Formula Based N-Body Algorithms},
  author={Hu, Y Charlie and Johnsson, S Lennart},
  institution={Center for Research in Computing Technology, Harvard University},
  year={1996}
}

@article{nabors:1994:fmm,
author = {Nabors, K. and Korsmeyer, F. and Leighton, F. and White, J.},
title = {Preconditioned, Adaptive, Multipole-Accelerated Iterative Methods for Three-Dimensional First-Kind Integral Equations of Potential Theory},
journal = {SIAM Journal on Scientific Computing},
volume = {15},
number = {3},
pages = {713-735},
year = {1994},
doi = {10.1137/0915046},
}

\end{document}
